%% file: afim-arxiv2.tex
\documentclass[12pt,a4paper]{article}

\usepackage{tikz,tikz-3dplot}
\usepackage{pgfplots}
\usetikzlibrary{external}                                               
\usetikzlibrary{3d}

\usetikzlibrary{plotmarks}

\usepackage{amsmath}    
\usepackage{amssymb}
\usepackage{graphicx}   
\usepackage{verbatim}   
\usepackage{color}      

\usepackage{diagbox}
\usepackage{smartref}

\usepackage{enumerate}
\usepackage{mathrsfs}
\usepackage{empheq}
\usepackage{bbold}

\usepackage[final]{showlabels}

\usepackage{amsthm}

\newcommand{\del}{\partial}
\renewcommand{\theta}{\vartheta}
\renewcommand{\phi}{\varphi}
\newcommand{\vecc}[2]{\left ( \begin{array}{c}#1\\#2\\ \end{array}\right )}

\newcommand{\dd}{\mathrm{d}}

\newcommand{\id}{\mathbb{1}}

\newcommand{\ii}{\mathbb{i}}

\newcommand{\advectionspeed}{u}
\newcommand{\cfl}{c}

\def \iconsize {0.4}
\newcommand{\stenOLO}{
	\begin{tikzpicture}[scale = \iconsize,very thick,  
		every circle node/.style={draw},
		every rectangle node/.style={draw},
		color = gray!50!white]
		\node[black,circle] (C) at (0,0) {};
		\node[] (L) at (-1,0) {};
		\node[] (R) at (1,0) {};
		\draw[] (L)--(C);
		\draw[] (R)--(C);
	\end{tikzpicture}
}

\newcommand{\stenLOO}{
	\begin{tikzpicture}[scale = \iconsize,very thick,  
		every circle node/.style={draw},
		every rectangle node/.style={draw},
		color = gray!50!white]
		\node[circle] (C) at (0,0) {};
		\node[black] (L) at (-1,0) {};
		\node[] (R) at (1,0) {};
		\draw[] (L)--(C);
		\draw[] (R)--(C);
	\end{tikzpicture}
}

\newcommand{\stenOOL}{
	\begin{tikzpicture}[scale = \iconsize,very thick,  
		every circle node/.style={draw},
		every rectangle node/.style={draw},
		color = gray!50!white]
		\node[circle] (C) at (0,0) {};
		\node[] (L) at (-1,0) {};
		\node[black] (R) at (1,0) {};
		\draw[] (L)--(C);
		\draw[] (R)--(C);
	\end{tikzpicture}
}

\newcommand{\stenLLO}{
	\begin{tikzpicture}[scale = \iconsize,very thick,  
		every circle node/.style={draw},
		every rectangle node/.style={draw},
		color = gray!50!white]
		\node[black,circle] (C) at (0,0) {};
		\node[black] (L) at (-1,0) {};
		\node[] (R) at (1,0) {};
		\draw[black] (L)--(C);
		\draw[] (R)--(C);
	\end{tikzpicture}
}

\newcommand{\stenOLL}{
	\begin{tikzpicture}[scale = \iconsize,very thick,  
		every circle node/.style={draw},
		every rectangle node/.style={draw},
		color = gray!50!white]
		\node[black,circle] (C) at (0,0) {};
		\node[] (L) at (-1,0) {};
		\node[black] (R) at (1,0) {};
		\draw[] (L)--(C);
		\draw[black] (R)--(C);
	\end{tikzpicture}
}

\newcommand{\stenLOL}{
	\begin{tikzpicture}[scale = \iconsize,very thick,  
		every circle node/.style={draw},
		every rectangle node/.style={draw},
		color = gray!50!white]
		\node[circle] (C) at (0,0) {};
		\node[black] (L) at (-1,0) {};
		\node[black] (R) at (1,0) {};
		\draw[] (L)--(C);
		\draw[] (R)--(C);
	\end{tikzpicture}
}

\newcommand{\stenLLL}{
	\begin{tikzpicture}[scale = \iconsize,very thick,  
		every circle node/.style={draw},
		every rectangle node/.style={draw},
		color = gray!50!white]
		\node[black,circle] (C) at (0,0) {};
		\node[black] (L) at (-1,0) {};
		\node[black] (R) at (1,0) {};
		\draw[black] (L)--(C);
		\draw[black] (R)--(C);
	\end{tikzpicture}
}

\newcommand{\exLOOimOLL}{{\tiny\begin{tabular}{c} \stenOLL\\\stenLOO \end{tabular}}}
\newcommand{\exOOLimOLL}{{\tiny\begin{tabular}{c} \stenOLL\\\stenOOL \end{tabular}}}
\newcommand{\exOOLimLOL}{{\tiny\begin{tabular}{c} \stenLOL\\\stenOOL \end{tabular}}}
\newcommand{\exOLOimLOL}{{\tiny\begin{tabular}{c} \stenLOL\\\stenOLO \end{tabular}}}
\newcommand{\exLOOimLOL}{{\tiny\begin{tabular}{c} \stenLOL\\\stenLOO \end{tabular}}}
\newcommand{\exOLOimOLL}{{\tiny\begin{tabular}{c} \stenOLL\\\stenOLO \end{tabular}}}
\newcommand{\exOLOimLLL}{{\tiny\begin{tabular}{c} \stenLLL\\\stenOLO \end{tabular}}}
\newcommand{\exLOOimLLL}{{\tiny\begin{tabular}{c} \stenLLL\\\stenLOO \end{tabular}}}
\newcommand{\exOOLimLLL}{{\tiny\begin{tabular}{c} \stenLLL\\\stenOOL \end{tabular}}}
\newcommand{\exLOLimLOL}{{\tiny\begin{tabular}{c} \stenLOL\\\stenLOL \end{tabular}}}
\newcommand{\exOLLimLLL}{{\tiny\begin{tabular}{c} \stenLLL\\\stenOLL \end{tabular}}}
\newcommand{\exLOLimLLL}{{\tiny\begin{tabular}{c} \stenLLL\\\stenLOL \end{tabular}}}
\newcommand{\exLLOimLLL}{{\tiny\begin{tabular}{c} \stenLLL\\\stenLLO \end{tabular}}}
\newcommand{\exOLOimLLO}{{\tiny\begin{tabular}{c} \stenLLO\\\stenOLO \end{tabular}}}
\newcommand{\exOOLimLLO}{{\tiny\begin{tabular}{c} \stenLLO\\\stenOOL \end{tabular}}}
\newcommand{\exLOOimLLO}{{\tiny\begin{tabular}{c} \stenLLO\\\stenLOO \end{tabular}}}

\tikzset{external/only named=true}

\renewcommand{\title}{Implicit Active Flux methods for linear advection}

\newcommand{\authorOne}{Wasilij Barsukow\footnote{Bordeaux Institute of Mathematics, Bordeaux University and CNRS/UMR5251, Talence, 33405 France, wasilij.barsukow@math.u-bordeaux.fr}}
\newcommand{\authorTwo}{Raul Borsche\footnote{University of Kaiserslautern-Landau, Gottlieb-Daimler-Straße 48, 67663 Kaiserslautern, Germany, borsche@mathematik.uni-kl.de}}

\begin{document}

\begin{center} \Large
\title

\vspace{1cm}

\date{}
\normalsize

\authorOne, \authorTwo
\end{center}

\begin{abstract}

In this work we develop implicit Active Flux schemes for the scalar advection equation. At every cell interface we approximate the solution by a polynomial in time. This allows to evolve the point values using characteristics and to update the cell averages using fluxes obtained by integrating this polynomial. The resulting schemes have order of convergence up to five, but show only moderate oscillations with high frequencies for discontinuous solutions. In numerical experiments we compare the different methods and show an application to network flows.

Keywords: linear advection, implicit methods, Active Flux

Mathematics Subject Classification (2010): 65M06, 65M25

\end{abstract}

\section{Introduction}

Linear advection is the simplest hyperbolic PDE and is widely used as a starting point for the development of numerical methods for conservation laws. 
It is the perfect testbed for studying properties of numerical methods, e.g. the analysis of linear (von Neumann) stability or the order of convergence.
Apart from being a prototype for nonlinear problems or systems of equations, there are some applications relying directly on the advection equation \cite{leveque01,GoettlichSupply,BorscheHeating} or the wave equation in 1-d, which can be diagonalized with characteristic variables.

A special class of these problems considers advection phenomena on networks.
Scalar equations are used for modeling supply chains \cite{GoettlichSupply} or district heating systems \cite{BorscheHeating}.
The wave equation on networks is considered e.g. in \cite{Zuazua,EggerWave,BorscheWave,LeugeringWave}.
All these applications demand highly accurate and efficient numerical methods and there is an overwhelming amount of possible schemes \cite{leveque01,ToroBook}.
However, most of these schemes are explicit and have to obey some kind of CFL condition bounding the size of the time step relative to the advection speed. In industrial applications, the relevant networks can consist of many hundred edges with the ratio between the longest and the shortest edge $10^{-3}$ and less. Consider, for example, the longest edge to have unit length, to be discretized with 100 cells and the time step chosen such that $\text{CFL} = 1$ on it. An edge of length $10^{-3}$ cannot be discretized by less than one cell, and thus with the same flow speed on it $\text{CFL} = 10$. Imposing $\text{CFL} = 1$ on a short edge decreases efficiency and accuracy on the entire network. The challenges therefore are similar to the small-cell-problem occurring ubiquitously in applications where the mesh is imposed externally. Observe that the difficulties are solved by having stable methods with only moderately large $\text{CFL}$ conditions (below 10, say).
In order to avoid the time step restriction an implicit scheme shall be chosen \cite{Helluy2019,Formaggia19,Xin2011,Baeder2007,Harten1981}.

In this paper we want to extend the ideas from standard finite difference/finite volume methods to the 
recently developed Active Flux method.
Active Flux was elaborated in \cite{eymann11,eymann13}, but it began its existence as a method for linear advection, when it was proposed in \cite{vanleer77} as ``Scheme V''.  
It has been extended to numerical methods for other conservation laws, e.g. the Euler equations of ideal hydrodynamics, see \cite{kerkmann18,barsukow19activeflux}.
Its main difference compared to classical finite differences or finite volume approaches is that it evolves point values and cell averages simultaneously. Besides having industrial applications in mind, the aim of the present work is to investigate strategies for the development of implicit Active Flux methods for linear problems. We believe that this paper is also the starting point for future investigations towards 
implicit Active Flux methods for more complex problems.

Active Flux shall be briefly reviewed next.

Consider a one-dimensional equidistant grid with cells $[x_{i-\frac12},x_{i+\frac12}]$, $i \in \mathbb Z$ and spacing $\Delta x$. 
The degrees of freedom of Active Flux are cell averages $\{\bar q_i \}_{i\in\mathbb Z}$ and point values $\{ q_{i+\frac12} \}_{i \in\mathbb Z}$ located at cell interfaces such that
\begin{align*}
	\bar q_i(t) \simeq \frac{1}{\Delta x} \int_{x_{i-\frac12}}^{x_{i+\frac12}} q(t, x) \,\dd x,
	\qquad &\qquad 
	q_{i+\frac12}(t) \simeq q(t, x_{i+\frac12}).
\end{align*}

The point values are evolved independently of the averages, contrary to e.g. the parabolic spline method \cite{zerroukat05}, where the point values are computed at each time step from the given averages. The evolution of the cell averages is immediately obvious: Integrating a conservation law $\del_t q + \del_x f(q) = 0$ over the cell yields
\begin{align}
 \frac{\dd}{\dd t}\bar q_i(t) + \frac{f(q_{i+\frac12}(t)) - f(q_{i-\frac12}(t))}{\Delta x} &= 0 \label{eq:averageupdatesemidiscrete}
\end{align}
The order of accuracy of the update of the average is entirely given by the order of accuracy of the point value update, as equation \eqref{eq:averageupdatesemidiscrete} is exact.

There exist several suggestions for the point value update, all of them being explicit in time:
\begin{enumerate}[(i)]
\item \label{it:rk} In \cite{abgrall20,abgrall22} it has been proposed to replace the space derivative by a (suitably upwinded) finite difference that uses the point values and averages, and thus to write down a semidiscretization of the conservation law for the point value. Together with \eqref{eq:averageupdatesemidiscrete}, the two equations can then be updated in time using a standard, explicit (e.g. Runge-Kutta) method. Depending on the choice of the finite difference, these methods are stable for CFL numbers well below 1.

\item \label{it:char} Initially, it was proposed in \cite{vanleer77} to define a parabolic reconstruction in every cell, whose average matches $\bar q_i$ and which passes through $q_{i\pm\frac12}$ at the endpoints of the cell, and to use it as initial data for a characteristics-based update of the point values. This means that the point value $q_{i+\frac12}(t^{n+1})$ at some time $t^{n+1}$ is found at the foot point of the characteristic which passes through $x_{i+\frac12}$ at time $t^{n+1}$. This ensures upwinding and stability. The natural CFL condition prevents the foot point of the characteristic from being further away than in a neighbouring cell. Von Neumann stability results agree with this ``physical'' stability bound, see \cite{vanleer77,chudzik21} for further details.
\end{enumerate}

Whatever the point value update, virtually all suggestions so far were time-explicit. 
In deriving implicit methods, we follow both strategies (\ref{it:rk}) and (\ref{it:char}) in order to be able to compare the results:
\begin{enumerate}[(i)]
\item As in \cite{abgrall20,abgrall22} we propose to integrate a semidiscretization of the conservation law for the point value, as well as \eqref{eq:averageupdatesemidiscrete} implicitly in time using standard methods.

\item We propose new implicit one-stage Active Flux methods. The derivation of these methods is, on the one hand, based on the usage of characteristics as in (\ref{it:char}) above, on the other hand on the idea of reconstruction in time, that was used for finite difference methods in \cite{MR4409827}. We analyze all methods in a certain class and identify those which are stable. 
\end{enumerate}
We find that the methods resulting from the latter approach are largely superior, which might be due to the fact that the time and space discretizations are not separated.

In \cite{nishikawa16}, Active Flux is used to solve a singular hyperbolic system that approximates the heat equation. For steady problems in this context, iterative solvers are used to solve the residual equation. For unsteady problems, it is suggested to incorporate a discrete time derivative into the residual. This means that ``a pseudo-steady problem [...] needs to be solved at every physical time step'' (p. 76 in \cite{nishikawa16}). This is prohibitive for purely hyperbolic problems that are addressed in the present work, and we aim at avoiding any sub-iterations. 
We also let ourselves inspire by the fact that the point value update used in Active Flux mixes time and space discretization. 
This is in contrast to standard methods where a space discretization is chosen first, and then a standard discretization in time is imposed. In principle, it could be that important properties are lost, because the two discretizations do not ``match'' (\cite{roe21}). In Section \ref{sec:numerical} we do indeed find that our proposed methods perform better than a method of lines.

Note that our methods have compact stencils due to the additional degrees of freedom of Active Flux.
This is particularly advantageous whenever boundary conditions are imposed. 
Active Flux offers the other great advantage of having a point value located just at the cell interface, where a Dirichlet boundary condition can be imposed immediately and unambiguously. 

The paper is organized as follows. Section \ref{sec:RK-schemes} presents implicit Active Flux methods based on strategy (\ref{it:rk}), and Section \ref{sec:onestepchar} presents those obtained through characteristics and reconstruction in time, following (\ref{it:char}). In Section \ref{sec:boundary} the implementation of Dirichlet boundaries is discussed, which demonstrates the advantages associated to using an Active Flux method. Numerical examples are shown in Section \ref{sec:numerical}.

\section{Implicit Active Flux schemes}

We aim at solving the scalar advection equation with fixed positive speed $\advectionspeed \in \mathbb R$, $\advectionspeed > 0$
\begin{align}\label{eq:advection}
	\partial_t q + u \partial_x q &=0, & q &\colon \mathbb R^+_0 \times I \to \mathbb R, 
\end{align}
on a compact domain $I \subset \mathbb R$, endowed with either periodic or an inflow Dirichlet boundary condition and an initial condition $q_0 \colon I \to \mathbb R$, $q(0, x) = q_0(x)$.

\subsection{Semi-discrete methods integrated implicitly}
\label{sec:RK-schemes}

Semi-discrete Active Flux methods have been first introduced in \cite{abgrall20}. The evolution equation of the cell average is trivially given by \eqref{eq:averageupdatesemidiscrete}. In order to obtain an evolution equation for the point value, a finite difference approximation to the spatial derivative is used. A third-order approximation is (see \cite{abgrall22})
\begin{align}
\label{eq:pointupdatesemidiscretethird}
 \frac{\dd}{\dd t} q_{i+\frac12}(t) &= - \advectionspeed \frac{2 q_{i-\frac12}(t) - 6 \bar q_i(t) + 4 q_{i+\frac12}(t)}{\Delta x}. \end{align}
Note that there is no notion of a conservative update for a point value, and the only condition that needs to be imposed is stability. Equations \eqref{eq:averageupdatesemidiscrete} and \eqref{eq:pointupdatesemidiscretethird} are a coupled system of ODEs that can be solved with standard methods. Explicit Runge-Kutta schemes were used in \cite{abgrall22}; here the system shall be integrated in time implicitly.

\newcommand{\sfrac}[2]{#1/#2}

We have studied the following methods:
\begin{itemize}
\item backward Euler
\item Crank-Nicolson 
\item Radau IA (3rd order) \quad $\displaystyle  \begin{array}{c|cc} 0 & \sfrac14 & - \sfrac14 \\ \sfrac23 & \sfrac14 & \sfrac5{12} \\ \hline & \sfrac14 & \sfrac34 \end{array}$

\item Radau IIA (3rd order) \quad $\displaystyle  \begin{array}{c|cc} \sfrac13 & \sfrac5{12} & - \sfrac1{12} \\ 1 & \sfrac34 & \sfrac14 \\ \hline & \sfrac34 & \sfrac14 \end{array} $

\item DIRK (Crouzeix, 3rd order) \quad $\displaystyle \begin{array}{c|cc} \sfrac12 + \sfrac{\sqrt{3}}6 & \sfrac12 + \sfrac{\sqrt{3}}6 & 0 \\ \sfrac12 - \sfrac{\sqrt{3}}6 & -\sfrac{\sqrt{3}}3 & \sfrac12 + \sfrac{\sqrt{3}}6 \\ \hline & \sfrac12 & \sfrac12 \end{array} $

\end{itemize}
for \eqref{eq:averageupdatesemidiscrete} and \eqref{eq:pointupdatesemidiscretethird}, and find all of them to be stable experimentally. 
Obviously only the latter methods yield the necessary (3rd) order of accuracy that corresponds to that of the spatial discretization. 
The Radau methods are of true multi-stage nature, which entails a significant increase in the linear system that needs to be solved at every time step. Numerical results are shown in Section \ref{sec:numerical}. 

\subsection{Single-stage methods}
\label{sec:onestepchar}

\subsubsection{Derivation}
\label{sec:derivation}

Classical explicit schemes are often based on a polynomial interpolation in space in each cell.
These polynomials are used to compute the fluxes at the interfaces, e.g. by solving a Riemann problem, which for linear advection amounts to tracing back the characteristics.
This procedure, however, does not pair well with large time steps, as the characteristics will move a distance larger than the size of a single cell. 
In principle, it is possible to keep track of the correct cell, which leads to large-time-step methods \cite{leveque85}. 

\begin{figure}[h]
\centering
\includegraphics[width=0.8\textwidth]{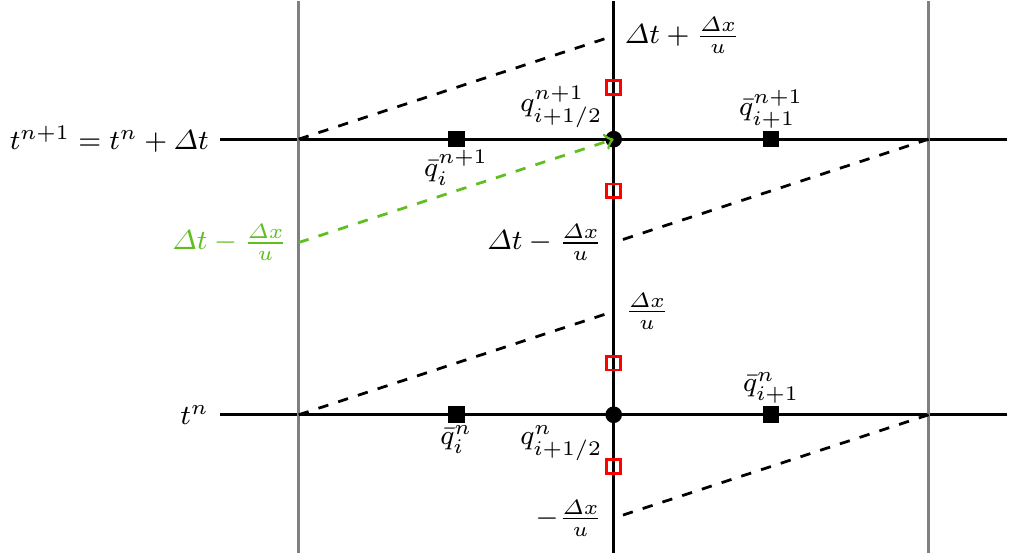}
\caption{Sketch of the construction of the polynomial at the interface. Characteristics are shown as dashed lines. Therefore, the average $\bar q_i^n$ (black square) is equal to the average $\frac{1}{\Delta t} \int_0^{\frac{\Delta x}{\advectionspeed}} q_{i+\frac12}^\text{recon}(t) \, \dd t$ (red square) of the reconstruction in time.}
	\label{fig:derivation}
\end{figure}

Here, however, to overcome this restriction we choose an interpolation in time at cell interfaces.
The values at different spatial locations are transported by the characteristics to the interface, which can involve values at both the current time level $t^n$ and also at the next time level $t^{n+1}$. 
These values are then interpolated polynomially, and the interpolation in time mixes values from both time levels into one polynomial at the interface.
The procedure is illustrated in Figure \ref{fig:derivation}.
In the following we give a detailed description of the construction of the methods.

We aim at obtaining an implicit numerical method that has a compact stencil.
Thus, for the interpolation at an interface we consider a combination of the values adjacent to the interface at both times $t^n$ and $t^{n+1}$.
For the interface at $x_{i+\frac{1}{2}}$ these are
\begin{align}
	 \bar q_i^n, q_{i+\frac12}^n,\bar q_{i+1}^n,\qquad 
	 \text{and}\qquad
	 \bar q_i^{n+1}, q_{i+\frac12}^{n+1},\bar q_{i+1}^{n+1}.
	 \label{eq:recondofs}
\end{align}
We therefore can choose among polynomials up to degree $5$ and thus construct schemes of at most $6^\text{th}$ order of accuracy.

Depending on the degree of the polynomial, we can choose among the following equations to find a reconstruction polynomial $q_{i+\frac12}^\text{recon}$ in time:
\begin{align}
 q_{i+\frac12}^\text{recon}(t^{n}) &= q_{i+\frac12}^{n} \label{eq:allequations1}\\
 q_{i+\frac12}^\text{recon}(t^{n+1}) &= q_{i+\frac12}^{n+1} \label{eq:allequations2} \\
 \frac{\advectionspeed}{\Delta x} \int^{t^{n+1} + \frac{\Delta x}{\advectionspeed}}_{t^{n+1}} q_{i+\frac12}^\text{recon}(t) \, \dd t &= \bar q_{i}^{n+1} 
\nonumber
 \\
 \frac{\advectionspeed}{\Delta x} \int^{t^{n+1}}_{t^{n+1} - \frac{\Delta x}{\advectionspeed}} q_{i+\frac12}^\text{recon}(t) \, \dd t &= \bar q_{i+1}^{n+1} \label{eq:allequations4}\\
 \frac{\advectionspeed}{\Delta x} \int^{t^{n} + \frac{\Delta x}{\advectionspeed}}_{t^n} q_{i+\frac12}^\text{recon}(t) \, \dd t &= \bar q_{i}^{n} 
 \nonumber
 \\
 \frac{\advectionspeed}{\Delta x} \int^{t^{n} }_{t^n - \frac{\Delta x}{\advectionspeed}} q_{i+\frac12}^\text{recon}(t) \, \dd t &=\bar  q_{i+1}^{n}. \label{eq:allequations6}
\end{align}
Note that the interpolation can exceed the time interval $[t^n,t^{n+1}]$, as indicated in Figure \ref{fig:derivation}.

Once the polynomial is determined, the update of the cell averages follows the classical update formula
\begin{align}
 \bar q_i^{n+1} = \bar q_i^n - \Delta t \frac{\hat f^{n+\frac12}_{i+\frac12} - \hat f^{n+\frac12}_{i-\frac12}}{\Delta x} \label{eq:finvolupdate2}
\end{align}
with 
\begin{align}
 \hat f^{n+\frac12}_{i+\frac12} := \frac{1}{\Delta t} \int^{t^{n+1}}_{t^n} f\left(q_{i+\frac12}^{\text{recon}}(t)\right) \, \dd t. \label{eq:fluxquadrfromrecon}
\end{align}
The point values can be updated by tracing back the characteristic 
\begin{align}
	q_{i+\frac32}^{n+1} = q_{i+\frac12}^\text{recon}\left(t^{n+1} - \frac{\Delta x}{\advectionspeed}\right). \label{eq:pointvalueupdatesinglestep}
\end{align}

Note that the reconstruction generally depends on values at the time level $t^{n+1}$. Thus, equations \eqref{eq:pointvalueupdatesinglestep} with \eqref{eq:fluxquadrfromrecon} and \eqref{eq:finvolupdate2} are implicit formulas, with the unknowns appearing linearly.

The stencil of the method involves at most the following values
\begin{align*}
 \bar q_i^n, q_{i+\frac12}^n\bar q_{i+1}^n,&&\bar q_i^{n+1}, q_{i+\frac12}^{n+1}\bar q_{i+1}^{n+1}, &&\text{ to update } q_{i+\frac32}^{n+1},\\
 \bar q_{i-1}^n, q_{i-\frac12}^n, \bar q_i^n, q_{i+\frac12}^n,\bar q_{i+1}^n ,&&
 \bar q_{i-1}^{n+1}, q_{i-\frac12}^{n+1}, \bar q_i^{n+1}, q_{i+\frac12}^{n+1},\bar q_{i+1}^{n+1}, &&
  \text{ to update } \bar q_{i}^{n+1}.
\end{align*}
Thus the update involves only values of neighboring cells. 

In order to easily refer to the different methods, we will use a \textbf{pictorial representation}. 
The 6 symbols (boxes/circles) in\exOLOimOLL represent the 6 degrees of freedom possibly involved in the reconstruction in time, as in \eqref{eq:recondofs}. The point values are represented by circles and the averages by boxes, with the upper row denoting time level $t^{n+1}$ (implicit) and the lower the time level $t^n$ (explicit). 
Finally, the black symbols are those actually in use for the reconstruction. The symbol shown involves a parabolic reconstruction fulfilling \eqref{eq:allequations1}, \eqref{eq:allequations2} and \eqref{eq:allequations4}.
For simplicity, in Section \ref{app:allmethods}, the stable methods are also given a unique \textbf{identifier} consisting of their order of accuracy and a capital letter.

Below we illustrate the construction of the schemes on a particular example.

\subsubsection{Example}

In this example we aim to construct a scheme of order three.
Thus the reconstruction polynomial $q_{i+\frac12}^\text{recon}$ has to be quadratic and we can 
choose 3 equations out of \eqref{eq:allequations1}--\eqref{eq:allequations6}. 
For example, one might select \eqref{eq:allequations1}, \eqref{eq:allequations2} and \eqref{eq:allequations4}, i.e. those involving $q_{i+\frac12}^{n}, q_{i+\frac12}^{n+1}, q_{i+1}^{n+1}$. 
Then, after applying the interpolation described in \ref{sec:derivation} and some calculations one finds
\begin{align*}
 q_{i+\frac12}^\text{recon}(t) &= q_{i+\frac12}^n + (t-t^n) \frac{ 2 \advectionspeed \left(q_{i+\frac12}^n-q_{i+\frac12}^{n+1}+3 \cfl \left(\bar q_{i+1}^{n+1} \cfl- q_{i+\frac12}^n+(1-\cfl)q_{i+\frac12}^{n+1}\right)\right)   }{ \cfl (3 \cfl-2) \Delta x } \\
 &+ (t-t^n)^2 \frac{ 3 \advectionspeed^2 \left(-2 \bar q_{i+1}^{n+1} \cfl+q_{i+\frac12}^n+(-1+2 \cfl) q_{i+\frac12}^{n+1}\right)}{\cfl (3 \cfl-2) \Delta x^2}
\end{align*}
where $\cfl = \frac{\advectionspeed \Delta t}{\Delta x}$ is the CFL number. From here, \eqref{eq:pointvalueupdatesinglestep} gives
\begin{align*}
 q_{i+\frac32}^{n+1} = \frac{6 \bar q_{i+1}^{n+1} (\cfl-1) \cfl+q_{i+\frac12}^n-\left(1-4 \cfl+3 \cfl^2\right) q_{i+\frac12}^{n+1}}{\cfl (3 \cfl-2)}
\end{align*}
Moreover, \eqref{eq:fluxquadrfromrecon} yields the numerical flux
\begin{align*}
 \hat f^{n+\frac12}_{i+\frac12} = \advectionspeed \frac{\bar q_{i+1}^{n+1} \cfl^2+(\cfl-1) \left(q_{i+\frac12}^n+q_{i+\frac12}^{n+1}-\cfl q_{i+\frac12}^{n+1}\right)}{3 \cfl-2}
\end{align*}
and inserting it into \eqref{eq:finvolupdate2} finally gives (having brought all the terms on one side of the equation)
\begin{align*}
 0 &= \bar q_i^n (2-3 \cfl)+\bar q_{i+1}^{n+1} \cfl^3-\bar q_i^{n+1} (2-3 \cfl+\cfl^3) \nonumber \\&+(\cfl-1) \cfl \Big (-q_{i-\frac12}^n+(\cfl-1) q_{i-\frac12}^{n+1}+q_{i+\frac12}^n+q_{i+\frac12}^{n+1}-\cfl q_{i+\frac12}^{n+1} \Big) 
\end{align*}
Following the notation introduced above, this method is denoted by \exOLOimOLL. 
The corresponding identifyer from Section \ref{app:allmethods} is 3\ref{method:exOLOimOLL}.

\subsubsection{Stability of single-stage methods}

Von Neumann, or $\ell^2$ stability analysis aims at quantifying whether Fourier modes $\exp(\ii k x)$ of spatial frequency $k \in \mathbb R$ are amplified or damped by the numerical method. To this end, the ansatz for the solution at continuous level is taken as
\begin{align}
 q(t, x) = \hat q(t) \exp(\ii k x) \label{eq:fouriercont}
\end{align}
and the ansatz for the numerical solution (on equidistant grids) as
\begin{align}
 Q_i^n := \vecc{q_{i+\frac12}^n}{\bar q_i^n} = \underbrace{\vecc{\hat q^n}{\hat {\bar q}^n}}_{=: \hat Q^n} \exp(\ii k i \Delta x) .
 \label{eq:fourierdisc}
\end{align}
Define $\beta := k \Delta x$.
Observe that despite dealing with merely a scalar equation, Active Flux has two kinds of degrees of freedom, which are evolved independently. The fact of having two distinct functions is mirrored by having two equations as well. On the one hand, the discrete Fourier transform \eqref{eq:fourierdisc} of the numerical method reduces \eqref{eq:finvolupdate2} and \eqref{eq:pointvalueupdatesinglestep} to the system 
\begin{align}
\hat Q^{n+1} = A^{-1}B \hat Q^{n}, 
\label{eq:onestepadvancedisc}
\end{align}
where $A \in \mathscr M^{2 \times 2}(\mathbb C)$ is associated with the implicit part, and $B \in \mathscr M^{2 \times 2}(\mathbb C)$ with the explicit one. The non-singularity of $A$ is a prerequisite of a solvable method, and shall be assumed. Note that $A$ and $B$ depend on $k$. The ansatz $\hat Q^{n} = \hat Q^0 z^n$ for some yet to be determined $z \in \mathbb C$ yields
\begin{align*}
 \hat Q^0 z = A^{-1}B \hat Q^0 
\end{align*}
i.e. $z$ must be an eigenvalue of $A^{-1}B$ (which depends on $k$).

On the other hand, considering the ansatz \eqref{eq:fouriercont} for the advection equation $\del_t q + \advectionspeed \del_x q = 0$ implies
\begin{align*}
 q(t+\Delta t, x) = \hat q(t + \Delta t) \exp(\ii k x) = q(t, x - \advectionspeed \Delta t) = \hat q(t) \exp(\ii k x) \exp(- \ii k \advectionspeed \Delta t)
\end{align*}
Thus, we have
\begin{align}
 \hat q(t + \Delta t) = \hat q(t) \exp(- \ii k \advectionspeed \Delta t) \label{eq:onestepadvancecont}
\end{align}
and consequently $|\hat q(t + \Delta t)| = |\hat q(t)|$. A natural stability requirement for the numerical methods therefore is $|z| \leq 1$ for both eigenvalues.

We have studied von Neumann stability of all the 20 methods of third, 15 methods of fourth, 6 methods of fifth and the unique method of sixth order that result from the procedure described above (a total of 42 methods). In many cases, the eigenvalues of the complex-valued $2 \times 2$ matrix $A^{-1}B$ could be determined analytically\footnote{One observes that the computation of the inverse is, in fact, unnecessary, since $0 = \det(A^{-1}B - z \id)$ is equivalent to $0 = \det A \det (A^{-1}B - z \id) = \det(B - A z)$.}. For example one finds for \exOOLimOLL (3\ref{method:exOOLimOLL}) the values $z = 0$ and
\begin{align*}
z &= - \frac{2+\cos\beta - \ii \cfl \sin \beta}{2-\cfl^2+\cos\beta+\cfl^2 \cos\beta+2 \ii \cfl \sin\beta}.
\end{align*}

In the remaining cases, we applied the algorithm of \cite{miller71} (originally due to Schur \cite{schur17,schur18}), which allows to determine whether the zeros of a polynomial are contained in the unit disc without actually computing them. We applied the algorithm to a sampling of values of $\beta$ and for $\cfl < 10$.

It is also possible to analyze the stability of the methods by considering a fixed grid size $N$ with periodic boundaries and analyzing the eigenvalues of the $2N \times 2N$ update matrices. We also performed this type of stability analysis for CFL numbers ranging from $1$ to $10$ on a grid of 100 cells. The two methods of stability analysis gave the same results.

Among the methods studied, i.e. among all the methods of the form \eqref{eq:finvolupdate2}--\eqref{eq:pointvalueupdatesinglestep}, whose reconstruction in time involves at most the degrees of freedom mentioned in \eqref{eq:recondofs}, 
\begin{itemize}
\item 1 is marginally stable,
\item 15 have some finite CFL number $\cfl_\text{min} \geq 0$ above which they are stable.
\end{itemize}

Among the 15 stable methods, 12 are stable for $\cfl >2$ or better, and 8 are stable for $\cfl > 1$ or better. 
One among them is even unconditionally stable. 
These results are summarized in Tables \ref{tab:cumulativestbilitycounts} and \ref{tab:stabilitysinglestepmethods}.

\begin{table}[h]
 \centering
 \begin{tabular}{c|c|c|c||c}
                            & order 3 & order 4 & order 5 & total \\ \hline
   stable $\forall \cfl$    &        &    1    &        &    1    \\\hline
   stable for $\cfl > 1$    &   4     &   3     &   1     &   8     \\\hline
   stable for $\cfl > 2$    &   6     &   3     &   3     &  12      \\\hline\hline
   stable above some $\cfl$ &   9     &    3    &    3    &  15      
 \end{tabular}
\caption{Numbers of stable methods by order of accuracy and type of CFL conditions. The divisions are cumulative (and thus not exclusive): for example a method with a minimum CFL number of 1 is also counted in ``stable for $\cfl > 2$'' and ``stable above some $\cfl$''.}
 \label{tab:cumulativestbilitycounts}
\end{table}

\newcommand{\unstable}{{\color{gray}$\times$}}
\setlength\tabcolsep{1.5pt}
\begin{table}[h]
	\centering
		\begin{tabular}{|c|c|c|c|c|c|c|c|}
			\hline
			\backslashbox{impl.}{expl.}&\raisebox{-0.3\height}{\stenLOO}&\raisebox{-0.3\height}{\stenOLO}&\raisebox{-0.3\height}{\stenOOL}&\raisebox{-0.3\height}{\stenLLO}&\raisebox{-0.3\height}{\stenLOL}&\raisebox{-0.3\height}{\stenOLL}&\raisebox{-0.3\height}{\stenLLL}
			\\	\hline
			\raisebox{-0.3\height}{\stenLOO}&  &  &  & \unstable & \unstable{\color{gray}$^*$} & \unstable{\color{gray}$^*$} & \unstable{\color{gray}$^*$} 
			\\	\hline
			\raisebox{-0.3\height}{\stenOLO}&  &  &  & \unstable{\color{gray}$^*$} & \unstable{\color{gray}$^*$} & \unstable{\color{gray}$^*$} & \unstable{\color{gray}$^*$} 
			\\	\hline
			\raisebox{-0.3\height}{\stenOOL}&  &  &  & \unstable{\color{gray}$^*$} & \unstable{\color{gray}$^*$} & \unstable{\color{gray}$^*$}  & \unstable{\color{gray}$^*$}  
			\\	\hline
			\raisebox{-0.3\height}{\stenLLO}& {$>4.74^*$} &$>4.55$&$>3.74^*$& \unstable{\color{gray}$^*$} & \unstable{\color{gray}$^*$} & \unstable{\color{gray}$^*$} & \unstable{\color{gray}$^*$} 
			\\	\hline
			\raisebox{-0.3\height}{\stenLOL}&$\notin [1,2]$&$>1$&$>1$& \unstable{\color{gray}$^*$} &$(\times)$& \unstable{\color{gray}$^*$} & \unstable{\color{gray}$^*$} 
			\\	\hline
			\raisebox{-0.3\height}{\stenOLL}&$>2$&$>1$&$>1$& \unstable & \unstable{\color{gray}$^*$} & \unstable{\color{gray}$^*$} & \unstable{\color{gray}$^*$} 
			\\	\hline
			\raisebox{-0.3\height}{\stenLLL}&$>1$&$>1.10$&$>0$&$>2^*$&$>2$&$>1$& \unstable 
			\\	\hline
		\end{tabular}
\caption[a]{CFL conditions for stable schemes (schemes marked with \unstable are unstable for sufficiently large $\cfl$). The pictorial notations indicate, which of the values $[\bar{q}_i,q_{i+1/2},\bar{q}_{i+1}]$ are used (black symbols) in the interpolation at the interface at $x_{i+1/2}$. The columns denote the choice of degrees of freedom at time $t^{n}$ (explicit), the rows denote the choice of degrees of freedom at time $t^{n+1}$ (implicit). For example, the method \exOLOimLLO (3\ref{method:exOLOimLLO}) is stable for $\text{CFL} > 4.55$. The methods marked by an asterisk $^*$ are (additionally) stable for some low CFL numbers. The method marked with $(\times)$ is marginally stable and found to be rather oscillative in practice. The methods left blank are of an order of accuracy less than 3, and were not studied.}
\label{tab:stabilitysinglestepmethods}
\end{table}

\subsubsection{Analysis of numerical diffusion and dispersion} \label{ssec:diffdisp}

So-called diffusion and dispersion errors of a linear numerical method convey the information how Fourier modes of different spatial frequencies $k$ are amplified/damped and how far their speed of propagation differs from the analytical value. By comparing \eqref{eq:onestepadvancecont} with \eqref{eq:onestepadvancedisc} one observes that the eigenvalues $z$ of $A^{-1}B$ (which are functions of $\beta := k \Delta x$ and $\cfl$) need to be compared to $\exp(- \ii k \advectionspeed \Delta t) = \exp(- \ii \beta \cfl)$. It thus makes sense to quantify the numerical diffusion by computing $|z|$, and the numerical dispersion by computing $\frac{\mathrm{arg }z}{-\beta \cfl}$, the analytic value being in both cases 1 for all values of $\beta$. 
This analysis is to quantify the behaviour of a numerical method beyond merely requiring stability $|z|  \leq 1$. 

\begin{figure}[h]
\centering
\includegraphics[width=0.8\textwidth]{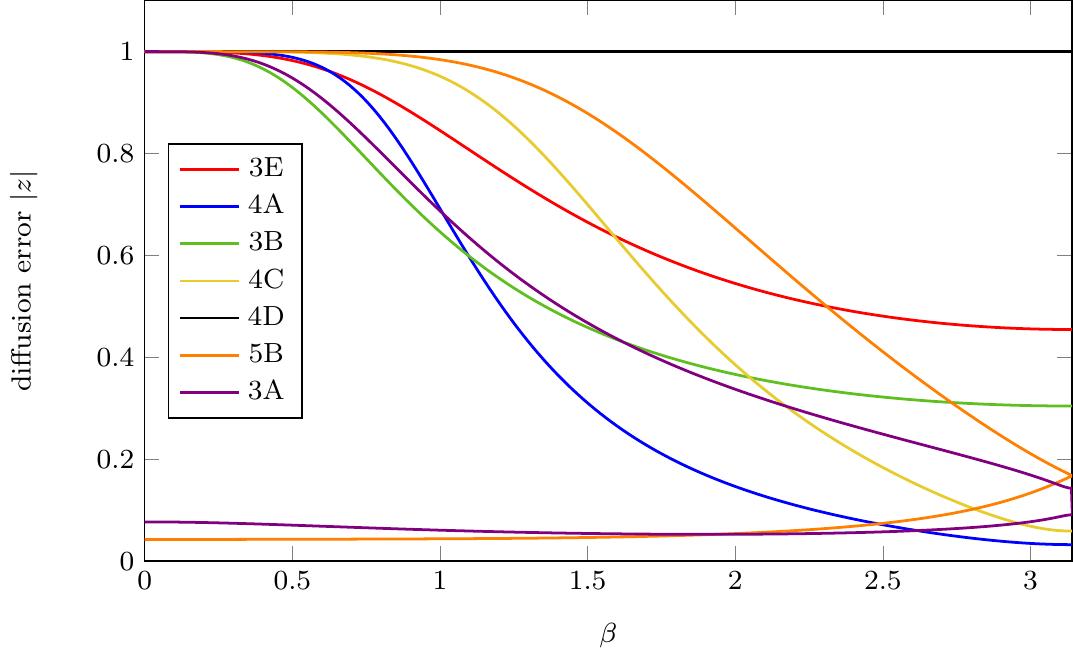}
	\caption[a]{Examples of diffusion curves $|z|$ for $\cfl = 3$, shown as functions of $\beta$ for different methods, using their identifiers from Section \ref{app:allmethods}. In principle, there are two eigenvalues per method, but often one of the eigenvalues is identically 0; it then is not shown here. Stability requires $|z| \leq 1$ for all $\beta$; the method \exLOLimLOL (4\ref{method:exLOLimLOL}) is marginally stable with the non-zero eigenvalue fulfilling $|z| = 1$ $\forall \beta$, i.e. equal to the analytical value. This, however, is not good, as this method is not damping waves moving at the wrong speed, i.e. those having large dispersion errors (see Figure \ref{fig:dispersion}).}
	\label{fig:diffusion}
\end{figure}

\begin{figure}[h]
\centering
\includegraphics[width=0.8\textwidth]{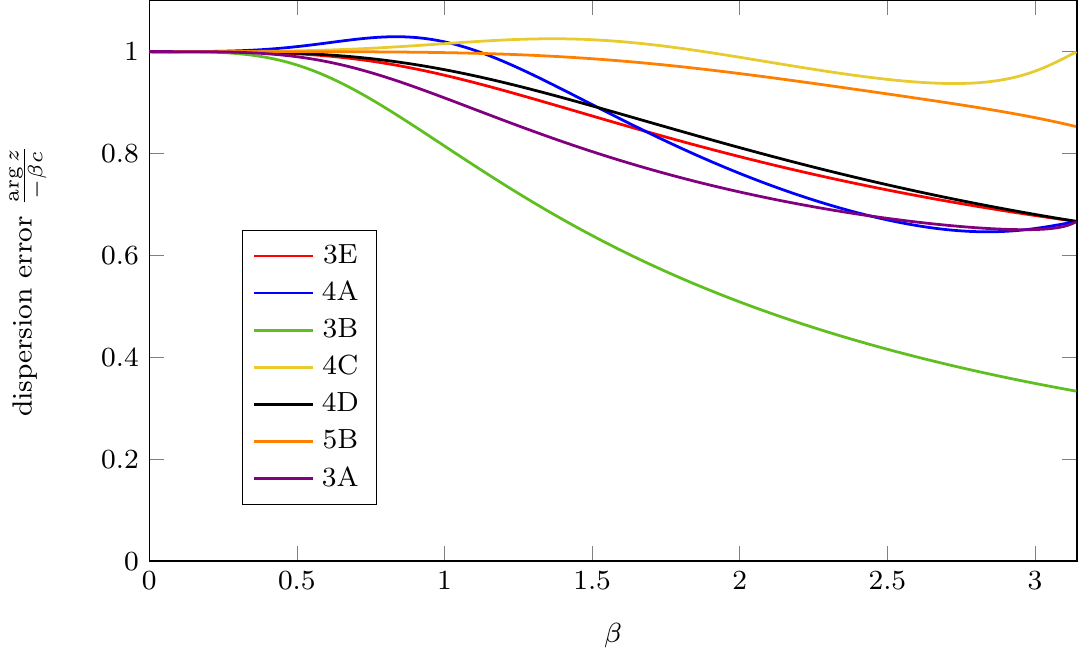}
	\caption[a]{Examples of dispersion curves $\frac{\arg z}{-\beta \cfl}$ for $\cfl = 3$, shown as functions of $\beta$ for different methods, using their identifiers from Section \ref{app:allmethods}. Only the physical eigenvalue is shown.}
	\label{fig:dispersion}
\end{figure}

Typical examples of diffusion and dispersion curves are shown in Figures \ref{fig:diffusion}--\ref{fig:dispersion}. Note that the eigenvalues depend on both $\beta$ and $\cfl$; the Figure shows $|z|$ and $\frac{\mathrm{arg }z}{-\beta \cfl}$ as functions of $\beta \in [0, \pi]$ for $\cfl = 3$. As there are two eigenvalues $z$, in principle, two curves appear for each method. One eigenvalue converges to the analytic one, whereas the other is spurious; it is sometimes zero, in which case it is not shown.

One observes that generally the error both in the diffusion and the dispersion increases with $\beta$. 
The diffusion is monotone, which is good, because then waves traveling at wrong speeds are damped. 
The importance of such behaviour has been emphasized in \cite{roe21} for explicit Active Flux methods and contrasted to that of other methods. 
The practical importance of damping those waves becomes obvious for the marginally stable method \exLOLimLOL (4\ref{method:exLOLimLOL}), which has $|z|=1$ for the non-vanishing eigenvalue and therefore does not damp waves which have wrong speeds. 
The corresponding numerical results (see Figure \ref{fig:Test_ord4}) show significantly more oscillations than the other stable methods.

In many cases a fine interplay of diffusion and dispersion is desirable.
Generally speaking, methods with less damping resolve sharp features better. 
Sharp features correspond to high values of $\beta$ and thus it is of interest to quantify how quickly the diffusion increases towards higher $\beta$. 
Recall that, given $\cfl$ and $\beta$, the time evolution of a Fourier mode is $|\hat Q^n| = |\hat Q^0| |z|^n = |\hat Q^0| |z|^{\frac{T}{\Delta t}}$. For example, consider $\advectionspeed = 1$, $\cfl = 3$, $\Delta x = 1/50$. Then $\Delta t = \frac{3}{50}$, and at $T = 8$ one has $|\hat Q^n| = |\hat Q^0| |z|^n = |\hat Q^0| |z|^{\frac{400}{3}}$. 
This is the setup of the numerical test in Section \ref{ssec:smoothanddiscprofiles}. 
The mode has decayed by half if $|z|^{\frac{400}{3}} = \frac12$, i.e. $|z|\simeq 0.995$. 

\begin{figure}[h]
\centering
\includegraphics[width=0.8\textwidth]{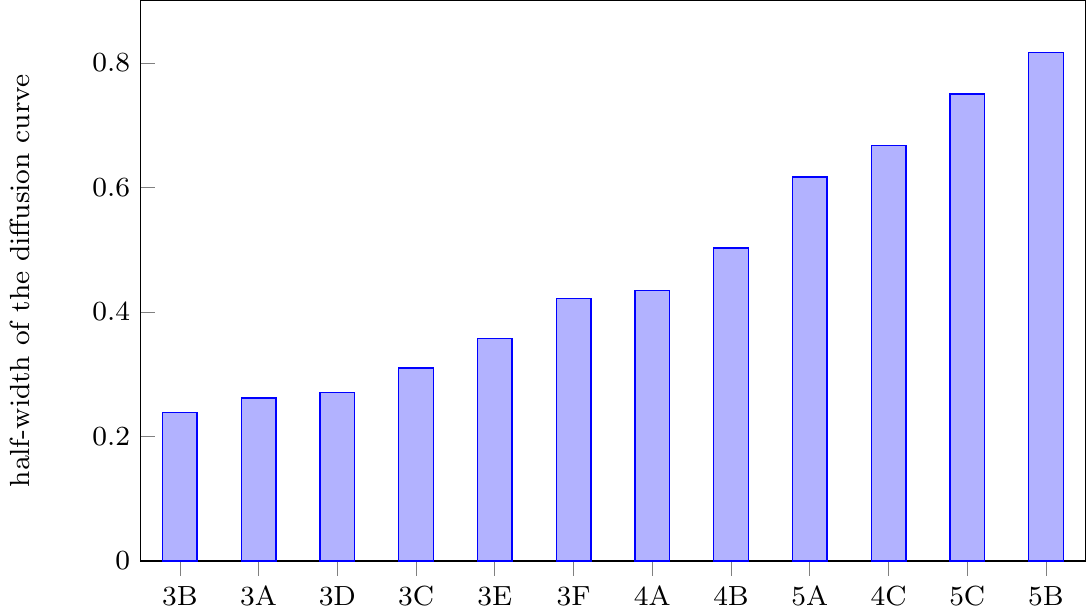}
	\caption[a]{Half-width of the diffusion curves at $|z| = 0.995$ for all the stable methods at $\cfl = 3$, i.e. the plot shows the value of $\beta$ for which $|z|$ attains the value $0.995$. The choice of values is explained in the text, and is made to fit the numerical tests shown in Section \ref{sec:numerical}. Thus, for values of $\beta$ less than the one shown the Fourier modes have not decayed by half by the time $T=8$. Identifiers from Section \ref{app:allmethods} are used to distinguish the different methods. The methods have been sorted in the order of increasing $\beta_{1/2}$, and thus methods further on the right are able to resolve sharp features better/for longer. Unsurprisingly, there is a correlation with the order of the method; however, for methods of the same order, large differences in this ability can be observed.}
	\label{fig:diffusionhalfwidth}
\end{figure}
It thus makes sense to compute, for a given method, the maximum ``frequency'' $\beta_{1/2}$ of the Fourier modes that has not yet decayed by half. 
The wavelength of these modes is $\frac{2 \pi \Delta x}{\beta_{1/2}}$. 
Measuring this frequency corresponds to measuring the half-width\footnote{The diffusion curves are bell-shaped and symmetric around $\beta = 0$. 
	The half-width is thus the distance along the abscissa between $\beta = 0$ and the location where the curve has a value $|z| = 0.995$.} of the diffusion curves (such as the ones shown in Fig. \ref{fig:diffusion}) at $|z| = 0.995$. 
These half-widths are shown in Figure \ref{fig:diffusionhalfwidth} for all the stable methods for $\cfl = 3$. From the plot one can see, for example, that for method \exOLOimOLL (3\ref{method:exOLOimOLL}) features on a length scale $\sim$20 $\Delta x$, corresponding to $\beta \simeq 0.31$, will have been diffused away by half by the time $T = 8$. 
For the method \exLOLimLLL (5\ref{method:exLOLimLLL}) the corresponding $\beta$ is 0.75 and the length scale $\sim 8 \Delta x$. This is in agreement with the numerical results of Section \ref{sec:numerical}.

We have also verified that for all but the marginally stable methods the eigenvalues of the amplification matrix decay to 0 as $\cfl \to \infty$ for all $k \neq 0$, which can be understood as L-stability of these methods.

\section{Boundary and coupling conditions}
\label{sec:boundary}

\subsection{Dirichlet boundary conditions}

On a finite interval $I = [x_\text{L}, x_\text{R}]$, the advection equation \eqref{eq:advection} with positive velocity $\advectionspeed>0$ has to be equipped with an initial condition $q_0$ and boundary data $b$ on the left end of the domain
\begin{align*}
	\begin{aligned}
	\partial_t q+\advectionspeed \partial_xq&=0\\
			q(t,x_\text{L})&=b(t) \quad \forall t > 0\\
			q(0,x)&=q_0(x) \quad \forall x \in I .
	\end{aligned}
\end{align*}
In the following we will discuss the modifications due to the boundary condition that need to be applied to the first two cells.

\begin{figure}[h]
\centering
\includegraphics[width=.45\textwidth]{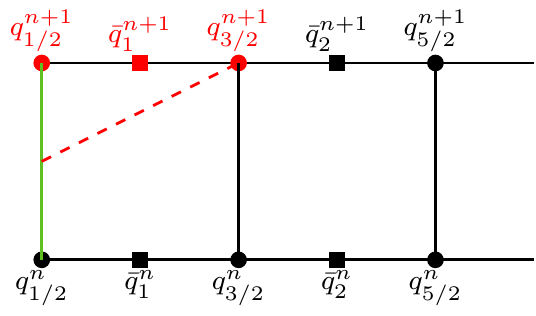}
	\caption{Scheme at the inflow boundary for the Dirichlet problem. The red values in the first cell can be read off directly by tracing back the characteristic (dashed lines) to the given boundary (green) at $x_{\frac12}$.}
	\label{fig:Dirichlet}
\end{figure}
\begin{figure}[h]
\centering
\includegraphics[width=.6\textwidth]{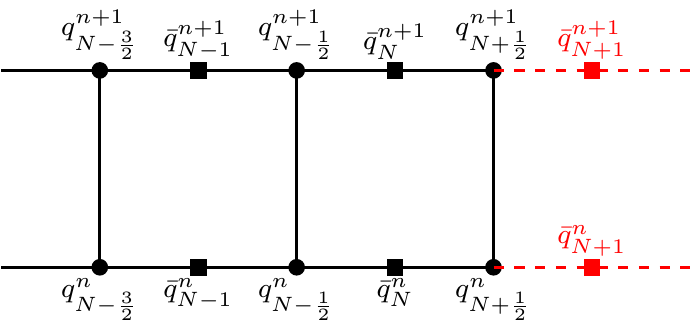}
	\caption{The outflow boundary. The values $\bar q_{N+1}^{n(+1)}$ (in red) are not available to define a reconstruction in time at $x_{N+\frac12}$. As explained in the text, in most cases it is not necessary to compute this reconstruction.}
	\label{fig:Dirichletoutflow}
\end{figure}

The main steps are illustrated in Figure \ref{fig:Dirichlet}.
The values $q^{n+1}_{\frac12}, q^{n+1}_{\frac32}$ of the first two interfaces, as well as the first cell average $\bar q^{n+1}_1$ are directly taken from the boundary data. 
By tracing back the characteristics we find
\begin{align*}
	q_{\frac12}^{n+1} = b(t^{n+1}),
	\quad 
	\bar{q}_{1}^{n+1} = \frac{\advectionspeed}{\Delta x}\int_{t^{n+1}-\frac{\Delta x}{\advectionspeed}}^{t^{n+1}}b(t)\dd t,
	\quad 
	q_{\frac32}^{n+1} = b\left(t^{n+1}-\frac{\Delta x}{\advectionspeed}\right)
	.
\end{align*}
Thus, it is easy to obtain the new values in the vicinity of the left (inflow) boundary by tracing the characteristics back to the boundary. 
The situation is different at the outflow end, where boundary conditions must not be given.
Assume that the grid is divided into $N$ cells (see Figure \ref{fig:Dirichletoutflow}). 
According to \eqref{eq:pointvalueupdatesinglestep} the right-most point value $q_{N+\frac12}$ is updated as
\begin{align}
 q_{N+\frac12}^{n+1} = q_{N-\frac12}^\text{recon}\left(t^{n+1} - \frac{\Delta x}{\advectionspeed}\right)
\end{align}
i.e. it only uses upwind information from the left because the reconstruction $q_{N-\frac12}^\text{recon}$ involves only the values $\bar q_{N-1}^{n(+1)}, q_{N-\frac12}^{n(+1)}, \bar q_{N}^{n(+1)}$ by definition. 

This is different for the average $\bar q^{n+1}_N$, because according to \eqref{eq:finvolupdate2} its update involves the quadrature of the two reconstructions in time $q^\text{recon}_{N-\frac12}$ and $q^\text{recon}_{N+\frac12}$. Three cases need to be distinguished:
\begin{enumerate}[1.]
 \item Assume that the reconstruction in time at $x_{i+\frac12}$ does not use the downwind cell averages $\bar{q}_{i+1}^{n(+1)}$. 
 Then, the reconstruction in time $q^\text{recon}_{N+\frac12}$ uses upwind values $\bar q_{N}^{n(+1)}, q_{N+\frac12}^{n(+1)}$ only and so does the update of the average $\bar q^{n+1}_N$. 
 Those schemes 
 can be applied directly. 
Further we note that such schemes can use an iterative procedure, marching the values from the left to the right: We first compute the update for the next point value and thereafter the update for the average.
This substitutes solving the linear system for the values at $t^{n+1}$, since the corresponding matrix is triangular.
 
 \item Assume next that the reconstruction in time at $x_{i+\frac12}$ involves the implicit downwind cell average $\bar{q}_{i+1}^{n+1}$ (and possibly its explicit counterpart), i.e. at $x_{N+\frac12}$ we need the non-available cell average $\bar{q}_{N+1}^{n+1}$. 
 The update equation for the last cell average according to \eqref{eq:finvolupdate2} reads
 \begin{align}
  \bar q_N^{n+1} = \bar q_N^n - \frac{\advectionspeed}{\Delta x}  \int_{t^n}^{t^{n+1}}   \left( q^\text{recon}_{N+\frac12}(t) - q^\text{recon}_{N-\frac12}(t)  \right)\, \dd t
   \label{eq:dirichletupdatelastequation}
 \end{align}
 where the right-hand-side depends, among other values, on $\bar q_{N+1}^{n+1}$. However, as the entire method is implicit, there is no notion of which equation updates which variable. Equation \eqref{eq:dirichletupdatelastequation} can equally well be seen as an update equation for $\bar q_{N+1}^{n+1}$ that involves only upwind values. What matters is that as many independent equations are provided as there are variables. 
 
 A different view of the same, inspired by the 'triangular'-scheme from the case above, is to say that the update equation for the cell average \eqref{eq:finvolupdate2} is not an update of $\bar{q}^{n+1}_{i}$ using the value $\bar{q}^{n+1}_{i+1}$, but an update of $\bar{q}^{n+1}_{i+1}$ using $\bar{q}^{n+1}_{i}$. Thus the matrix becomes triangular and all values can be computed in an iterative fashion.
 We thus propose to shift Equation \eqref{eq:dirichletupdatelastequation}, and in fact all the equations for the cell average, by one cell to the left. The index-shifted Equation \eqref{eq:dirichletupdatelastequation} would be counted as the equation updating $q_N^{n+1}$, and, in general, Equation \eqref{eq:finvolupdate2} would be counted as the update of the cell average $q_{i-1}^{n+1}$.

 The problem thus is shifted to the inflow boundary, where missing values are readily available by the procedure described earlier. For example, one can directly compute the integral of the flux
	$
	f_{1/2}^{n+1/2} = \int_{t^n}^{t^{n+1}}b(t)\dd t 
	$
as indicated in green in Figure \ref{fig:Dirichlet}.
 Which values are affected by this special treatment depends on the stencil of the reconstruction in time, but 
 any of them can be found in the Dirichlet boundary data upon tracing the characteristics. 
 Also, due to the compact stencil only a few degrees of freedom are updated directly using the boundary.
 
 \item At last we consider schemes using the explicit downwind cell average $\bar{q}_{i+1}^{n}$, but not the implicit one. The update equation \eqref{eq:dirichletupdatelastequation} would involve the unavailable cell average $\bar q_{N+1}^n$, but not $\bar q_{N+1}^{n+1}$. 
 It cannot thus, by shifting the index, be considered an update equation for $\bar q_{N}^{n+1}$ and the trick from above cannot be applied.
Fortunately, upon inspection of Table \ref{tab:stabilitysinglestepmethods}, one observes that there is only one stable (on domains with periodic boundaries) method which involves $\bar q_{i+1}^n$, but does not involve $\bar q_{i+1}^{n+1}$ in its reconstruction in time at $x_{i+\frac12}$: it is \exOOLimLLO (3\ref{method:exOOLimLLO}), and we skip searching for a fix.
\end{enumerate}

\subsection{Coupling conditions on networks}

In the case of networks, several edges are connected at a junction. 
Here the inflow boundary condition of an edge is determined by coupling conditions. Define the coordinates along an edge $e$ with length $\ell_e$ to be $x \in [0, \ell_e]$ with the orientation chosen such that the velocity points towards larger values of $x$.
Following \cite{MR4409827}, the inflow into an edge $e$ is determined by a weighted sum of the outflows from the edges at a junction $J$, i.e.
\begin{align*}
	q_{e}(t, 0)&=\alpha_{e}\sum_{e' \text{ flows into }J}q_{e'}(t, \ell_{e'})\qquad \text{for all }e\text{ with flow away from }J
\end{align*}
with 
\begin{align}
 \sum_{e\text{ flows away from }J} \alpha_{e} = 1 \qquad \text{for all }J,
\end{align}
where $\alpha_e \in [0,1]$ is the coupling weight of edge $e$, i.e. the proportion of mass available at the junction that actually enters $e$.

This can be interpreted as Dirichlet-type boundary conditions $b_e(t)$ when an interpolation in time $q_{e',N+\frac12}^\text{recon}(t)$ of the edges $e'$ which flow into the junction is known.
Denote by $\advectionspeed_{e'}$ and $\cfl_{e'}$ the speed of advection and the corresponding CFL number on edge $e'$. For third order of accuracy one can choose to reconstruct using $q_{e',N-\frac12}^n, q_{e',N+\frac12}^n, q_{e',N+\frac12}^{n+1}$, i.e.
\begin{align*}
 q_{e',N+\frac12}^\text{recon}(t^n) &= q_{e',N+\frac12}^n, & q_{e',N+\frac12}^\text{recon}(t^n + \Delta t) &= q_{e',N+\frac12}^{n+1}, \\
 q_{e',N+\frac12}^\text{recon}\left(t^n + \frac{\Delta x_{e'}}{\advectionspeed_{e'}}\right) &= q_{e',N-\frac12}^n.
\end{align*}
which gives
 \begin{align*}
  q_{e',N+\frac12}^\text{recon}(t) &= q_{e',N+\frac12}^n + \frac{t-t^n}{\Delta t} \frac{\cfl_{e'}^2 (q_{e',N+\frac12}^n - q_{e',N-\frac12}^n) + q_{e',N+\frac12}^{n+1} - q_{e',N+\frac12}^n}{1-\cfl_{e'}} \\
  \nonumber& + \frac{(t-t^n)^2}{\Delta t^2} \cfl_{e'} \frac{q_{e',N+\frac12}^{n+1} -  q_{e',N+\frac12}^n +\cfl_{e'} (q_{e',N+\frac12}^n - q_{e',N-\frac12}^n) }{\cfl_{e'} - 1}.
 \end{align*}
Choices involving other degrees of freedom are also possible, as well as additional degrees of freedom for higher order of accuracy. For fourth order, for example, we additionally use the average
\begin{align*}
\frac{\advectionspeed_{e'}}{\Delta x_{e'}} \int_{t^n}^{t^n + \frac{\Delta x_{e'}}{\advectionspeed_{e'}}} q_{e',N+\frac12}^\text{recon}\left(t\right) \,\dd t &= \bar q_{e',N}^n.
\end{align*}

Then, the Dirichlet boundary $b_e(t)$ on edge $e$ is
\begin{align*}
  b_e(t) = \,\,\,\,\alpha_e\!\!\!\!\! \sum_{e' \text{outflowing into }e} \!\!\!\!\! q_{e',N+\frac12}^\text{recon}(t).
\end{align*}

One then proceeds as indicated before, e.g.
\begin{align*}
 q_{e,\frac32}^{n+1} = b_e\left(\Delta t - \frac{\Delta x_e}{\advectionspeed_e} \right) .
\end{align*}

\section{Numerical tests} \label{sec:numerical}

\subsection{Convergence studies}

\begin{figure}[h]
	\centering
	\includegraphics[width=.8\textwidth]{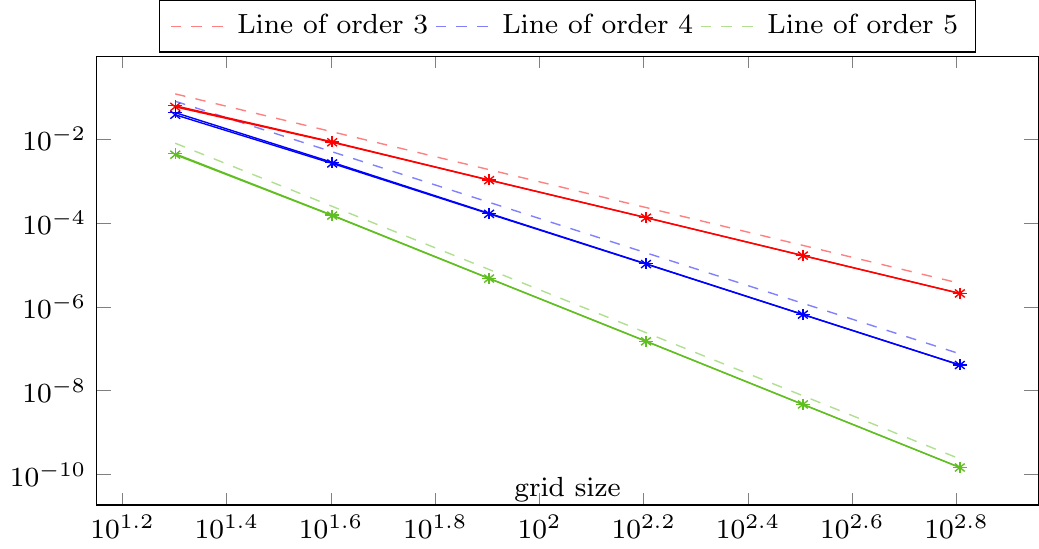}
		\caption[a]{Convergence tests. The $\ell^1$-error of the point values ($\times$) and the $L^1$-error of the averages ($+$) coincide optically with the dashed references. Shown are the schemes \\
		\exLOOimOLL (3\ref{method:exLOOimOLL}, \includegraphics{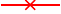}, 3$^\text{rd}$), \exOLOimLLL (4\ref{method:exOLOimLLL}, \includegraphics{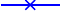}, 4$^\text{th}$) and \exOLLimLLL (5\ref{method:exOLLimLLL}, \includegraphics{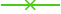}, 5$^\text{th}$).
	}
	\label{fig:convergence}
\end{figure}

We tested the order of convergence for all stable schemes discussed above.
The test with the initial condition $\sin(2\pi x)$ was run with a CFL number of about $3$ up to $T=10$.
The grids start with $20$ cells and are refined up to $640$ cells.
The errors on a grid with $N$ cells are computed as
\begin{align*}
	L^1=\Delta x \sum_{i=1}^{N}| \bar{q}_i - \bar q_i^\text{exact}| ,
	&&
	\ell^1 = \frac{1}{N}\sum_{i=1}^{N+1}|q_{i-\frac12} - q_{i-\frac12}^\text{exact}|.
\end{align*}
The results for some selected schemes are displayed in Figure \ref{fig:convergence}.
All schemes show the expected order of convergence and the $\ell^1$- and $L^1$-errors are always of similar size.

\subsection{Smooth and discontinuous profiles} \label{ssec:smoothanddiscprofiles}

Next we compare all the proposed methods in the same numerical test proposed in \cite{jiang96}, which combines smooth and discontinuous profiles.
For the Active Flux schemes we use a discretization with $100$ cells for the interval $[0,2]$ using periodic boundaries and a speed $\advectionspeed = 1$. 
This coincides with the number of degrees of freedom as proposed for the classical method on $200$ cells, which is what was used in \cite{jiang96}.
We run all methods with a CFL number of about $3$ up to a final time $T=8$, i.e. 4 grid revolution periods.

For reference we show the exact solution and the Finite Difference schemes of order $3$ and $4$ described in \cite{MR4409827,steinle93}.
Their solutions are computed on a grid with $200$ cells with the same time step as the Active Flux scheme (CFL number of about $6$), such that the memory and computational time are comparable. The 4th order Finite Difference method is marginally stable, which explains the large amount of oscillations.

\subsubsection{Semi-discrete method integrated implicitly}
\label{ssec:numericalradau}

\begin{figure}
\centering
\includegraphics[width=.8\textwidth]{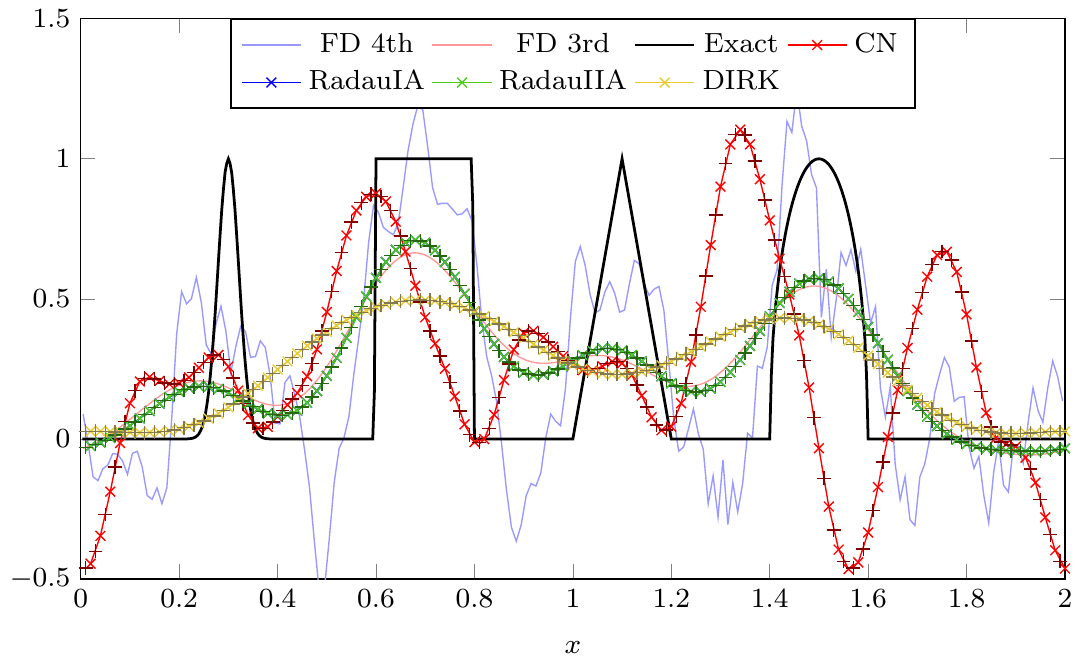}
	\caption{3rd order semi-discrete scheme integrated with common implicit methods. The curves corresponding to the two Radau methods are on top of each other.}
	\label{fig:Test_RK}
\end{figure}

Figure \ref{fig:Test_RK} shows results obtained by integrating the semi-discretization from section \ref{sec:RK-schemes} implicitly with standard methods.
With the second-order Crank-Nicholson method almost no features of the solution are preserved. As the spatial discretization of the semi-discretization \eqref{eq:pointupdatesemidiscretethird} is third-order accurate, only third-order time-integrators are of relevance.
The DIRK scheme is more diffusive than the reference third-order finite difference method, although third order itself.
Only the RadauIA and RadauIIA, which optically coincide, give results similar to the third order finite difference method.

\subsubsection{Third-order methods}
\begin{figure}
   \centering
   \includegraphics[width=.8\textwidth]{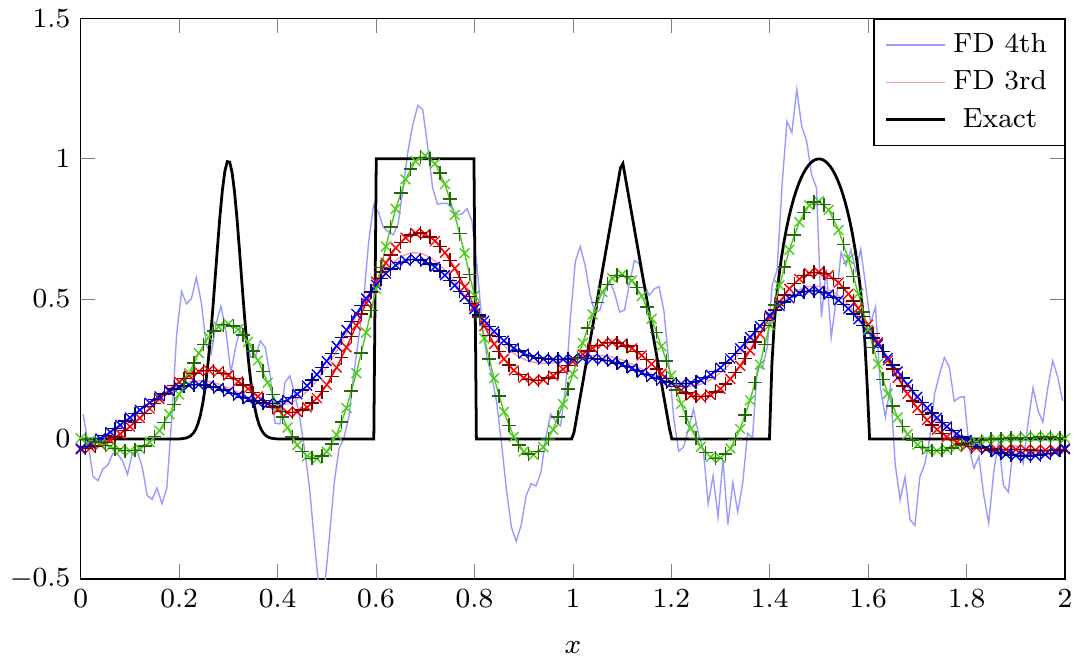}
	\caption[a]{Numerical test of the three third-order schemes ($\times$ point values, $+$ averages) \\ \exOOLimOLL (3\ref{method:exOOLimOLL},\includegraphics{tikz/ImplicitActiveFlux-figure_crossref_red.pdf}), \exOOLimLOL (3\ref{method:exOOLimLOL},  \includegraphics{tikz/ImplicitActiveFlux-figure_crossref_blue.pdf}) and \exLOOimOLL (3\ref{method:exLOOimOLL},  \includegraphics{tikz/ImplicitActiveFlux-figure_crossref_green.pdf}). 
		The third-order finite-difference method yields results very close to 3\ref{method:exOOLimLOL}.}
	\label{fig:Test_ord3_1}
\end{figure}

\begin{figure}
   \centering
   \includegraphics[width=.8\textwidth]{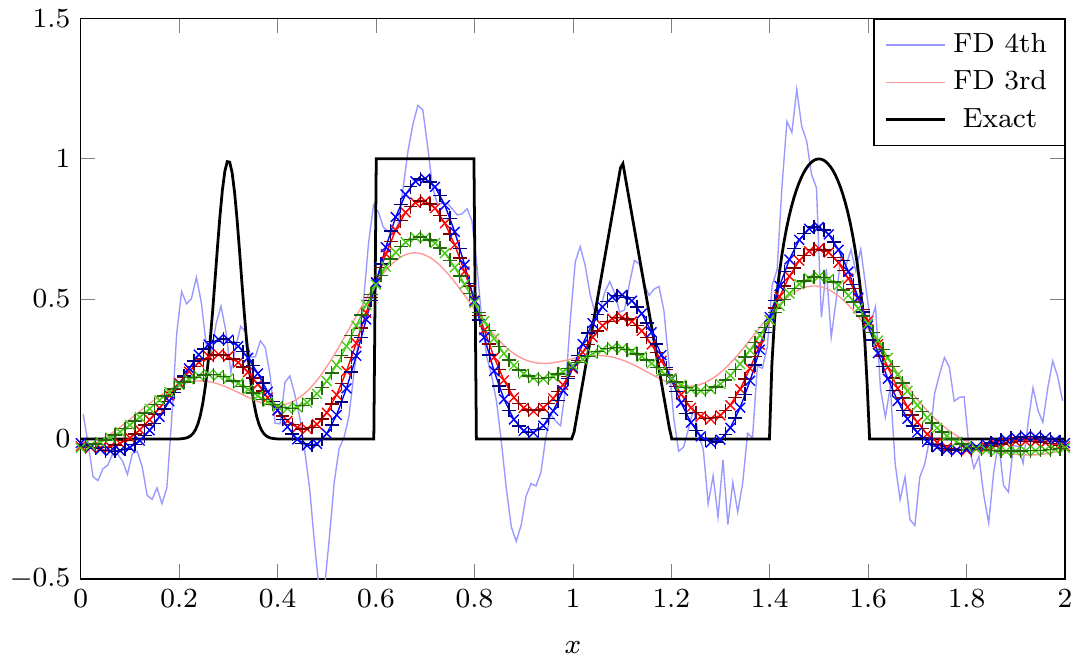}
	\caption[a]{Numerical test of the three third-order schemes ($\times$ point values, $+$ averages) \\
		\exOLOimOLL (3\ref{method:exOLOimOLL},  \includegraphics{tikz/ImplicitActiveFlux-figure_crossref_red.pdf}), \exLOOimLOL (3\ref{method:exLOOimLOL}, \includegraphics{tikz/ImplicitActiveFlux-figure_crossref_blue.pdf}) and \exOLOimLOL (3\ref{method:exOLOimLOL}, \includegraphics{tikz/ImplicitActiveFlux-figure_crossref_green.pdf}).}
	\label{fig:Test_ord3_2}
\end{figure}

\begin{figure}
   \centering
   \includegraphics[width=.8\textwidth]{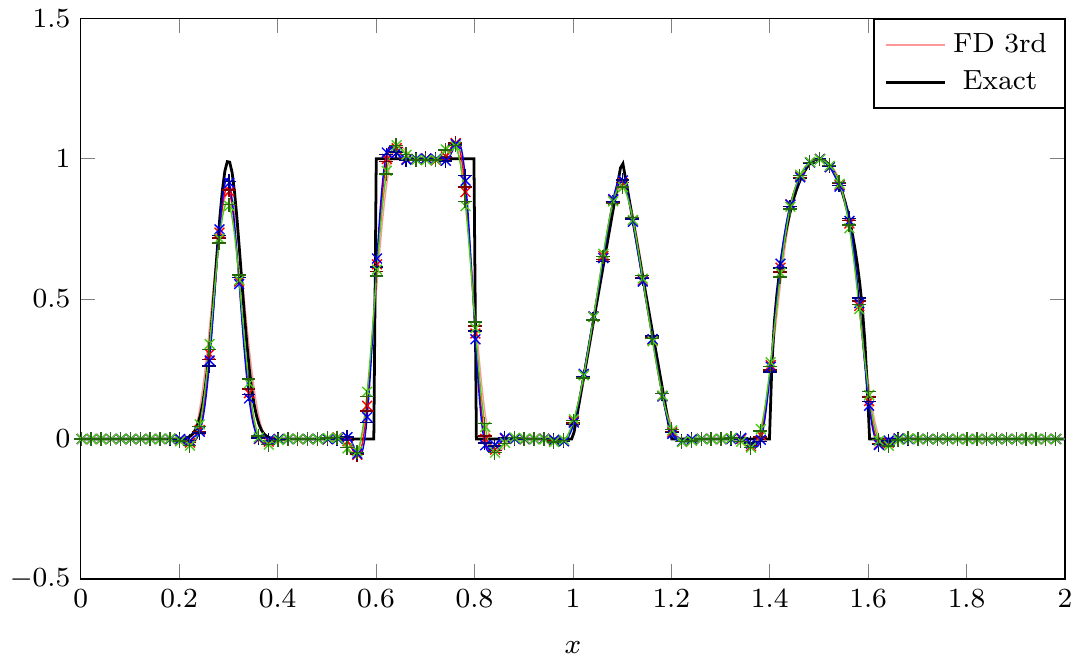}
	\caption[a]{Same as Figure \ref{fig:Test_ord3_2}, but on a grid of 1000 cells and marker every ten points.}
	\label{fig:Test_ord3_2_highres}
\end{figure}

In Figure \ref{fig:Test_ord3_1} and \ref{fig:Test_ord3_2} we show the schemes which use three values for the approximation at the interface and are thus of order $3$. We do not show \exOOLimLLO (3\ref{method:exOOLimLLO}), \exOLOimLLO (3\ref{method:exOLOimLLO}) and \exLOOimLLO (3\ref{method:exLOOimLLO}) (as they are unstable for CFL = 3).
All schemes perform better or equally well as the third order finite difference method.
The scheme \exOOLimLOL (3\ref{method:exOOLimLOL}) almost coincides with this reference.
All schemes show significant influence of numerical diffusion, but, mostly, the four solution features can be recognized.
Especially the methods with an upwind bias in the explicit part, i.e. \exLOOimOLL (3\ref{method:exLOOimOLL}) and \exLOOimLOL (3\ref{method:exLOOimLOL}), are significantly better compared to the other ones. They also clearly outperform the semi-discrete methods used with standard time-integrators (discussed in Section \ref{ssec:numericalradau}), especially noticing that two sub-stages of the Radau methods are required compared to a direct update. 
The differences between the results of the direct methods are less for larger CFL numbers (e.g. $\textrm{CFL}= 6$). Figure \ref{fig:Test_ord3_2_highres} shows the same setup upon refinement of the grid.

\subsubsection{Fourth-order methods} \label{ssec:numericalfourth}

\begin{figure}
   \centering
   \includegraphics[width=.8\textwidth]{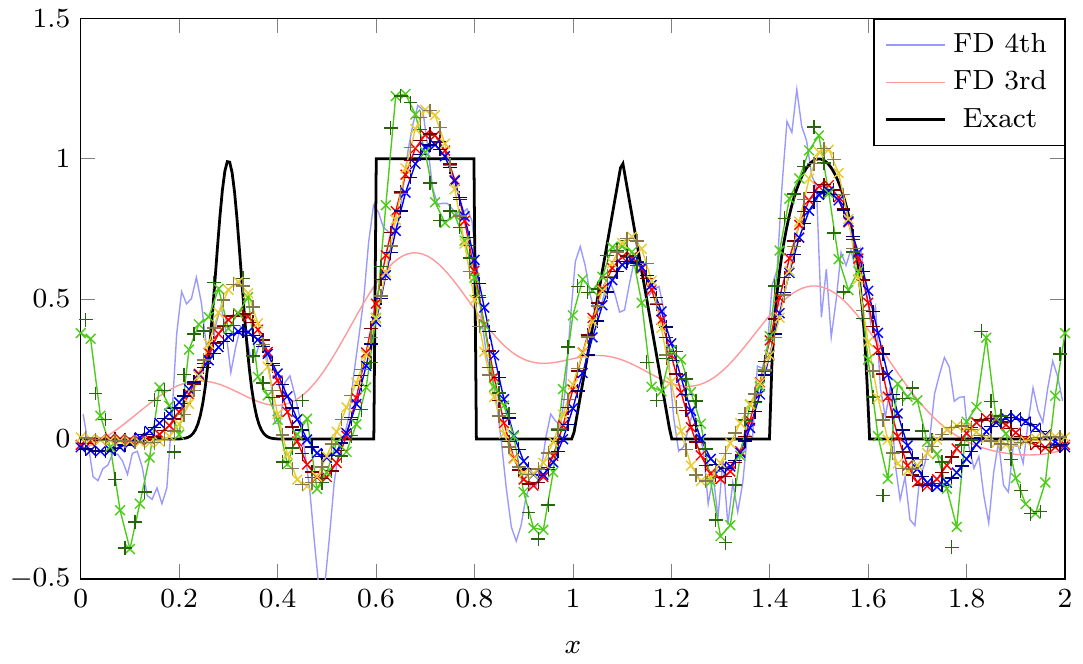}
	\caption[a]{Numerical test of the four fourth-order schemes ($\times$ point values, $+$ averages)\\ \exOLOimLLL (4\ref{method:exOLOimLLL}, \includegraphics{tikz/ImplicitActiveFlux-figure_crossref_red.pdf}), \exOOLimLLL (4\ref{method:exOOLimLLL}, \includegraphics{tikz/ImplicitActiveFlux-figure_crossref_blue.pdf}), \exLOLimLOL (4\ref{method:exLOLimLOL}, \includegraphics{tikz/ImplicitActiveFlux-figure_crossref_green.pdf}) and \exLOOimLLL (4\ref{method:exLOOimLLL}, \includegraphics{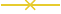}).}
	 \label{fig:Test_ord4}
\end{figure}

\begin{figure}
\centering
\includegraphics[width=.8\textwidth]{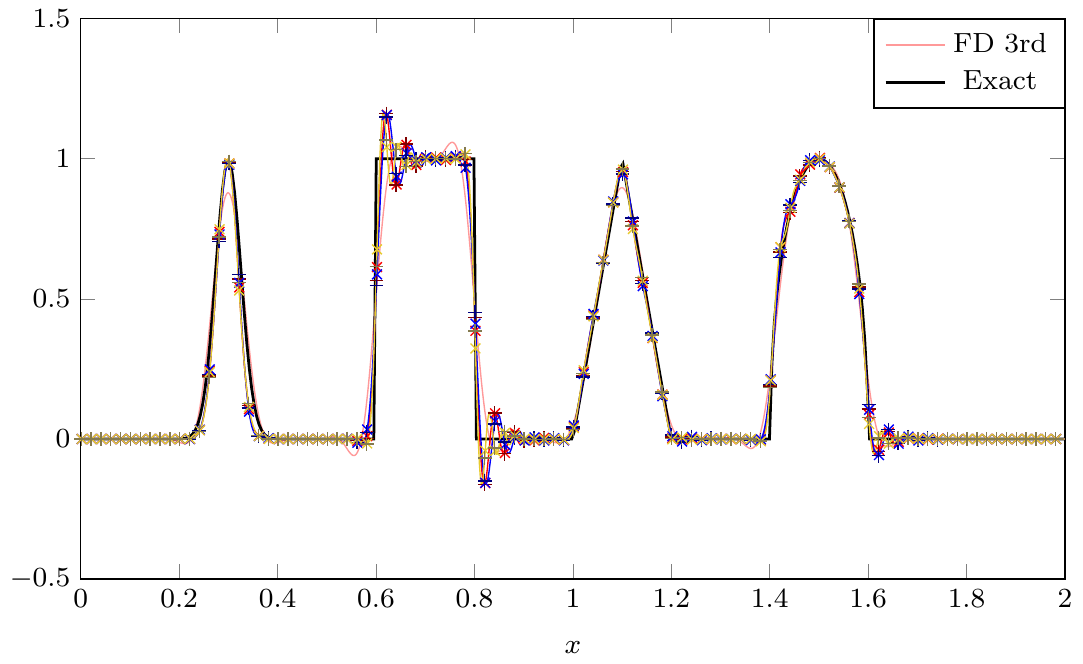}
	\caption[a]{Same as Figure \ref{fig:Test_ord4} but on a grid of 1000 cells and marker every ten points; without \exLOLimLOL (4\ref{method:exLOLimLOL}).}
	 \label{fig:Test_ord4_highres}
\end{figure}

In Figure \ref{fig:Test_ord4} four schemes of order 4 are shown. On coarse grids, all give more accurate results than third-order methods.
Scheme \exLOLimLOL (4\ref{method:exLOLimLOL}) plays a particular role, displaying significantly more oscillations than the other Active Flux methods. 
As can be seen in Table \ref{tab:stabilitysinglestepmethods}, it is marginally stable (the non-vanishing eigenvalue has norm 1), which might explain the oscillations. 
Moreover, the stencil of this method is such that the interpolation at the interface uses only average values and thus the update is independent of the point values.
It is therefore not surprising that method 4\ref{method:exLOLimLOL} behaves similarly to the fourth-order reference Finite Difference method.  
In the solutions of the other three methods the features can be identified much better.
Scheme \exLOOimLLL (4\ref{method:exLOOimLLL}), which has an upwind-bias in the explicit stencil seems to have a tendency of overshooting the exact solution. The main error of even-order methods is dispersive, and indeed upon grid refinement (Figure \ref{fig:Test_ord4_highres}) the methods display prominent high-frequency oscillations.

\subsubsection{Fifth-order methods}

\begin{figure}[h!]
   \centering
   \includegraphics[width=.8\textwidth]{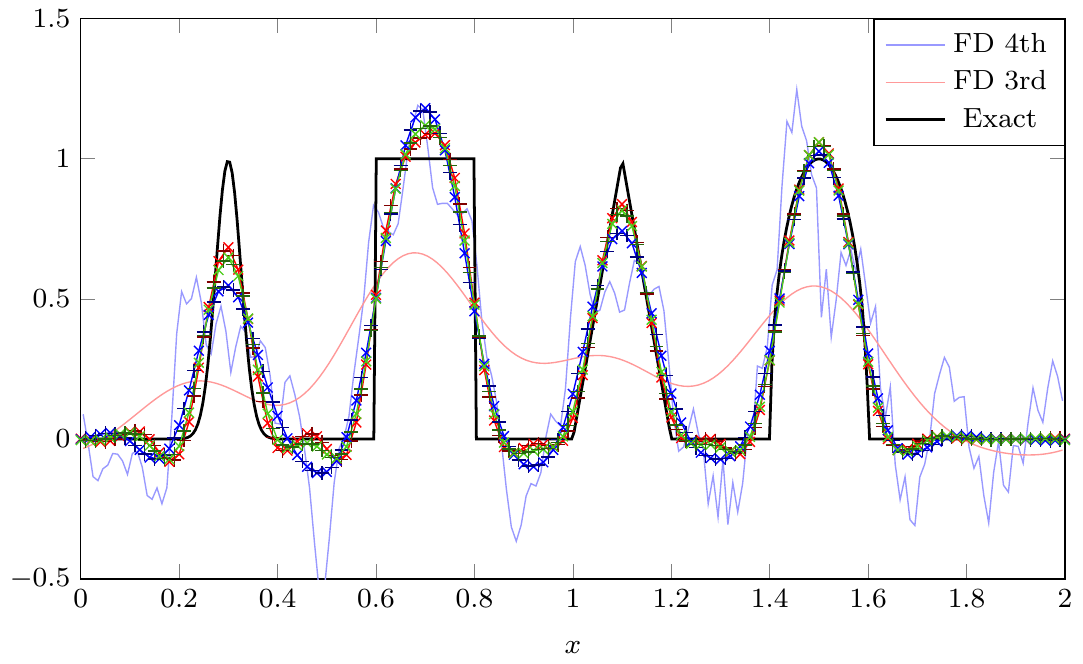}
	\caption[a]{Numerical test of the three fifth-order schemes ($\times$ point values, $+$ averages)\\ 
	\exLLOimLLL (5\ref{method:exLLOimLLL}, \includegraphics{tikz/ImplicitActiveFlux-figure_crossref_red.pdf}), \exOLLimLLL (5\ref{method:exOLLimLLL}, \includegraphics{tikz/ImplicitActiveFlux-figure_crossref_blue.pdf}) and \exLOLimLLL (5\ref{method:exLOLimLLL}, \includegraphics{tikz/ImplicitActiveFlux-figure_crossref_green.pdf}).}
	 \label{fig:Test_ord5}
\end{figure}

\begin{figure}[h!]
\centering
\includegraphics[width=.8\textwidth]{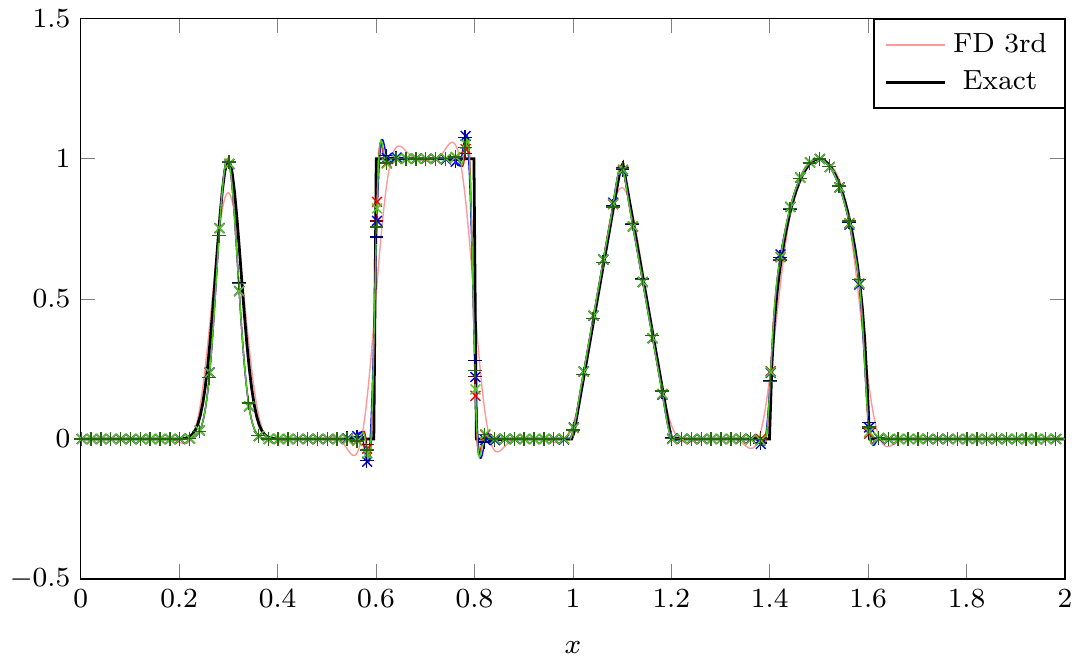}
	\caption[a]{Same as Figure \ref{fig:Test_ord5} but on a grid of 1000 cells and marker every ten points.}
	 \label{fig:Test_ord5_highres}
 \end{figure}

Finally, schemes of order 5 are shown in Figure \ref{fig:Test_ord5}.
One observes again a significant improvement in comparison to the fourth-order schemes.
All features of the solution are captured and can be clearly distinguished. The three schemes only differ in the choice of the explicit stencil, and a more upwind-focused stencil in the explicit part seems to improve the result slightly. Figure \ref{fig:Test_ord5_highres} shows results on a refined grid; as expected for odd-order methods the oscillations are moderate.

\subsubsection{Comparison to the explicit Active Flux method}

\begin{figure}[h!]
   \centering
   \includegraphics[width=.8\textwidth]{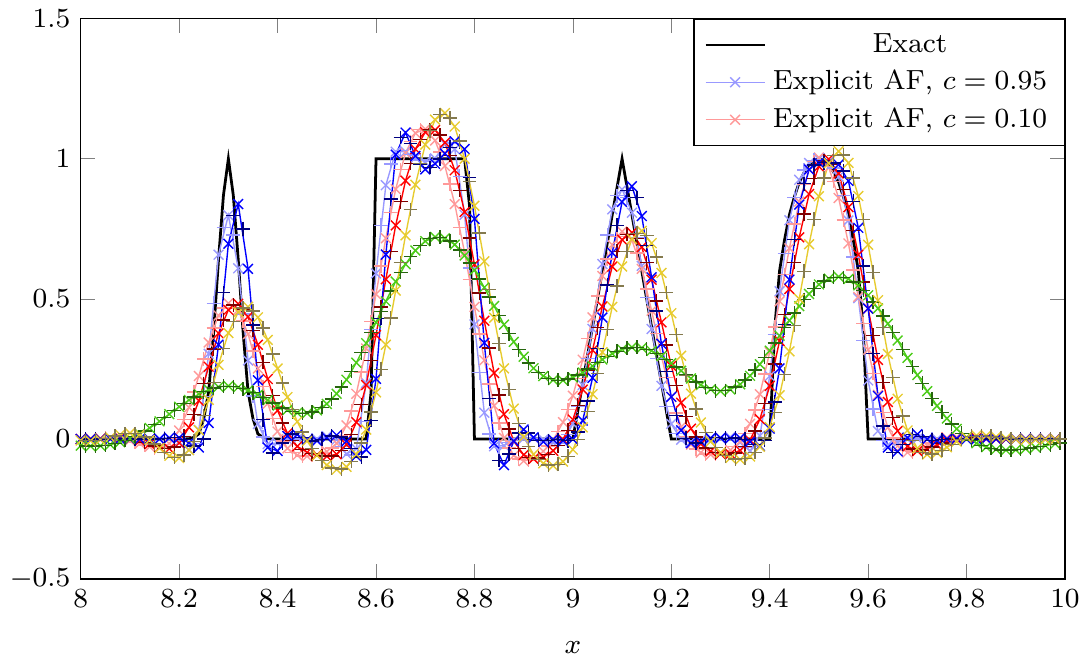}
	\caption[a]{Same setup as the Figure before, but now using an extended domain (keeping same $\Delta x$ and final time) and Dirichlet boundaries. This allows to implement the implicit methods with as much computational cost (per time step) as implicit ones. Here, the following methods are compared with the third-order explicit Active Flux method ($\times$ point values, $+$ averages):
		\exOLOimLOL (3\ref{method:exOLOimLOL}) $\cfl = 1.2$: \includegraphics{tikz/ImplicitActiveFlux-figure_crossref_red.pdf}, $\cfl = 3$: \includegraphics{tikz/ImplicitActiveFlux-figure_crossref_green.pdf}; \exOLLimLLL (5\ref{method:exOLLimLLL}) $\cfl = 1.2$: \includegraphics{tikz/ImplicitActiveFlux-figure_crossref_blue.pdf}, $\cfl = 3$: \includegraphics{tikz/ImplicitActiveFlux-figure_crossref_yellow.pdf}.}
	 \label{fig:comparisonexplicit}
\end{figure}

Figure \ref{fig:comparisonexplicit} shows a comparison between explicit Active Flux and some of the proposed implicit methods for different CFL numbers. As is well-known, numerical error vanishes for $\cfl = 1$. This is unachievable for practical applications. Still, one might be aiming at being close to $\cfl=1$ wherever one is interested to resolve details of the flow without much numerical diffusion. One deduces from Figure \ref{fig:comparisonexplicit} that our new implicit Method \exOLLimLLL (5\ref{method:exOLLimLLL}) with $\cfl = 1.2$ achieves similarly accurate results as (third-order) explicit Active Flux with a CFL of 0.95. The runtimes are also similar due to marching in space (0.06 sec for the former, and 0.05 sec for the latter) -- despite the implicit method being 5th order accurate. An increase in CFL number also increases the diffusion (with runtimes around 0.04 sec for the other implicit methods). 

The real advantage of using these implicit methods is related to them solving the short-edge-problem that occurs in the simulation of flow on a network, as it was briefly sketched in the Introduction. On ``interesting'' edges the CFL will be chosen close to 1 without much difference in computational effort between explicit and implicit methods. However, on short edges, this value will amount to, say, $\cfl=10$. Admittedly, there the solution will be smeared out, but this is acceptable for a short edge that itself is much smaller than computational cells on other edges. Also, the alternative is much more problematic: Upon usage of an explicit method one would be forced to choose $\cfl=0.1$ on the ``interesting'' edges in order to fulfill the $\cfl < 1$ condition on the short edge. This small CFL number will lead to both increased diffusion everywhere in the network and increased computational time. The increased diffusion of the explicit Active Flux method with $\cfl =0.1$ is visible in Figure \ref{fig:comparisonexplicit}, and the computation took 0.23 sec.

\subsection{Network}

\begin{figure}[h!]
\centering
\includegraphics[width=0.4\textwidth]{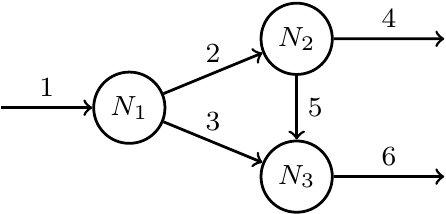}
 \caption{Sketch of the network with three nodes and six edges.}
 \label{fig:network}
\end{figure}

As a final test we consider a network of six edges and four nodes, shown in Figure \ref{fig:network}.
The lengths of the edges are
\begin{align*}
	\ell_1 = 5,\quad \ell_2=\ell_3=\ell_5 = 20,\quad \ell_4=\ell_6 = 30.
\end{align*}
On each edge we solve the advection equation with the speeds
\begin{align*}
	\advectionspeed_1 = \advectionspeed_3 = \advectionspeed_4 = \advectionspeed_6 = 1,
	\quad \advectionspeed_2 = 2,\quad \advectionspeed_5 = \frac{20}{10+\frac32} = \frac{40}{23}, 
\end{align*}
oriented according to the direction of the edges.\\
At the nodes $N_1,N_2,N_3$ suitable coupling conditions have to be imposed.
At the two splitting nodes $N_1$ and $N_2$ we distribute the incoming flux according to fixed parameters $\alpha_1 = \frac{3}{4}$ and $\alpha_2 = \frac{2}{3}$ such that
\begin{align*}
	q_2(t, 0) &= \alpha_1 q_1(t, \ell_1), & q_3(t, 0) &= (1-\alpha_1) q_1(t, \ell_1),\\
	q_4(t, 0) &= \alpha_2 q_2(t, \ell_1), & q_5(t, 0) &= (1-\alpha_2) q_2(t, \ell_1).
\end{align*}
The coupling at node $N_3$ is immediate by the conservation of mass
\begin{align*}
	q_6(t, 0) = q_5(t, \ell_5) + q_3(t, \ell_3).
\end{align*}
Finally, we impose at the first edge the Dirichlet boundary condition
\begin{align*}
	q_1(t, 0) = b(t) = \sin\left(\Omega t\right)
\end{align*}
with $\Omega = \frac{2\pi}{3}$ and an initial condition
\begin{align*}
	q_1(0, x) &= \exp\left(-4(x - \ell_1/2)^2\right)    & q_e(0, x) &= 0 \quad \mathrm{ for } \quad e=2, \ldots, 6
\end{align*}

Note that the configuration is chosen such that a signal starting in edge 1 will be split at node $N_1$ into two parts.
One signal is traveling directly along edge $3$ to node $N_3$, while the other takes a detour via $N_2$.
Apart from losing some of its strength to edge $4$, it arrives delayed at node $N_3$ as compared to the signal from edge $3$.

By defining the edge-crossing times $\tau_e := \ell_e / \advectionspeed_e$, $e = 1, \ldots, 6$, the exact solution for large times $t$ at the nodes $N_2$ and $N_3$ has the form
\begin{align*}
	q_4(t,0) &= \alpha_1 \alpha_2 b(t - \tau_1 - \tau_2),\\
	q_6(t,0) &= \alpha_1 (1-\alpha_2) b(t - \tau_1 - \tau_2 - \tau_5) + (1-\alpha_1) b(t - \tau_1 - \tau_3).
\end{align*}
Since the parameters $\alpha_1,\alpha_2$ and $\Omega$ are chosen such that they satisfy
\begin{align*}
	\alpha_1 (1-\alpha_2) &= 1-\alpha_1
\\
	\Omega (\tau_2 + \tau_5) &= \Omega \tau_3 + \pi 
\end{align*}
we have a destructive interference at $q_6(t,0)$.
This means that the exact solution on edge 6 is equal to zero once the pulse from the initial condition has passed.

\begin{figure}
\centering
\includegraphics[width=0.49\textwidth]{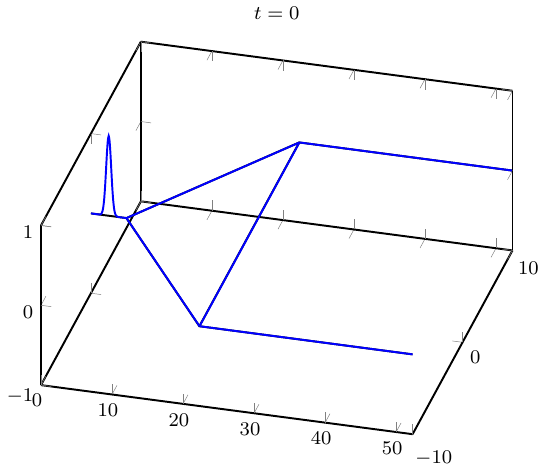} \hfill
\includegraphics[width=0.49\textwidth]{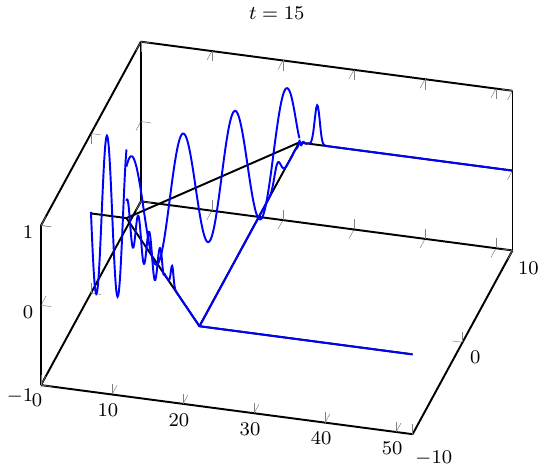}
\includegraphics[width=0.49\textwidth]{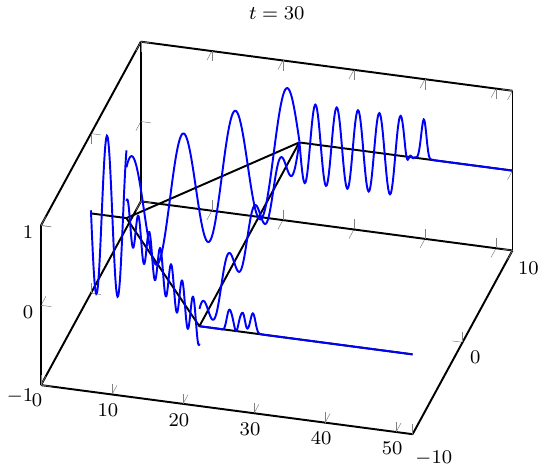} \hfill
\includegraphics[width=0.49\textwidth]{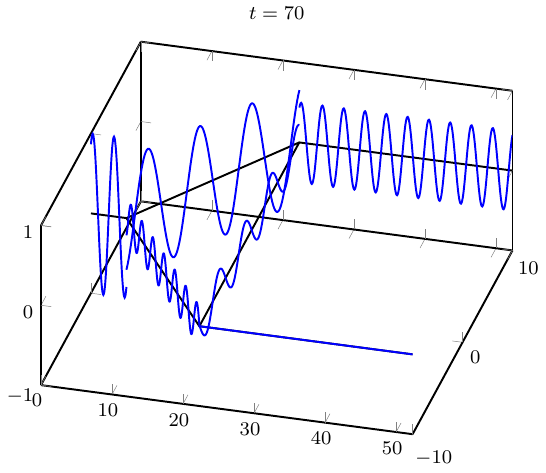}
 \caption[a]{Numerical solution of the flow on the network using Method \exOLOimLLL (4\ref{method:exOLOimLLL}) with $\text{CFL} = 5$ and $\Delta x = \frac18$ on all edges. 
	}
 \label{fig:networksimu}
\end{figure}

Figure \ref{fig:networksimu} shows the results of a simulation using Method \exOLOimLLL (4\ref{method:exOLOimLLL}) with $\text{CFL} = 5$ and $\Delta x = \frac18$ on all edges. 

\begin{figure}
\centering
\includegraphics[width=.8\textwidth]{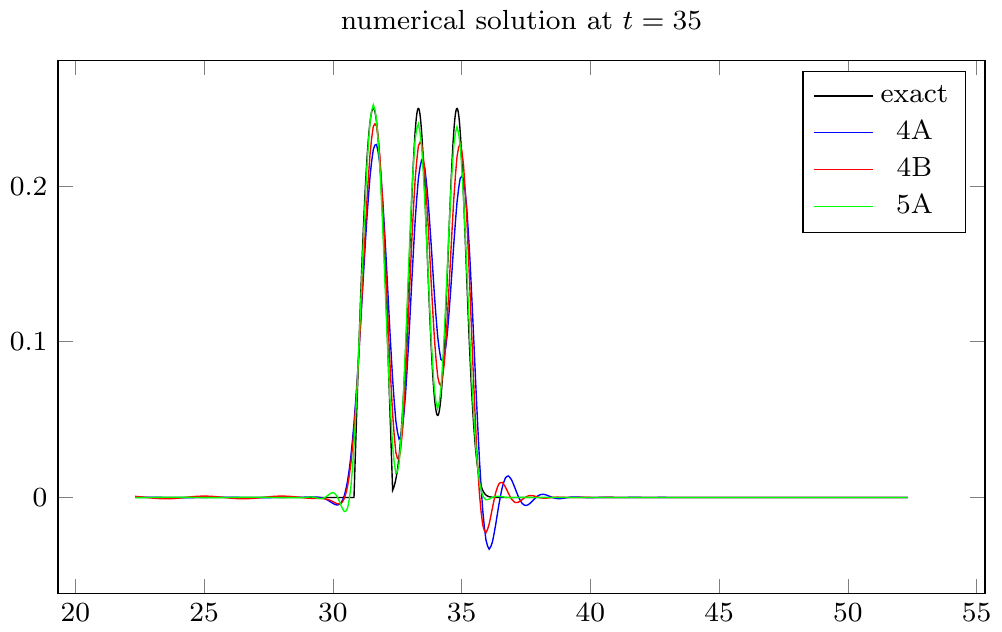}
	\caption[a]{Numerical solution on edge 6 at $t=35$ for different methods ($\Delta x = \frac18$).}
	\label{fig:networkconvhalftime}
\end{figure}

\begin{figure}
\centering
\includegraphics[width=.8\textwidth]{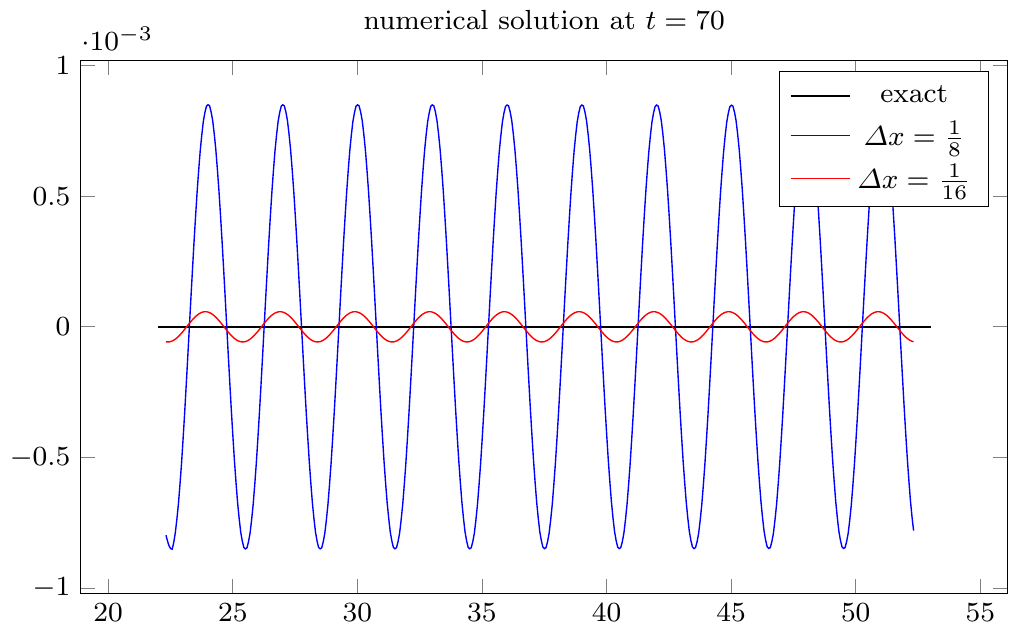}
	\caption[a]{Numerical solution on edge 6 at $t=70$. The exact solution is vanishing due to destructive interference. One observes the convergence of the numerical solution as $\Delta x$ is decreased from $\frac18$ to $\frac1{16}$. }
	\label{fig:networkconv}
\end{figure}

\begin{figure}
\centering
\includegraphics[width=.8\textwidth]{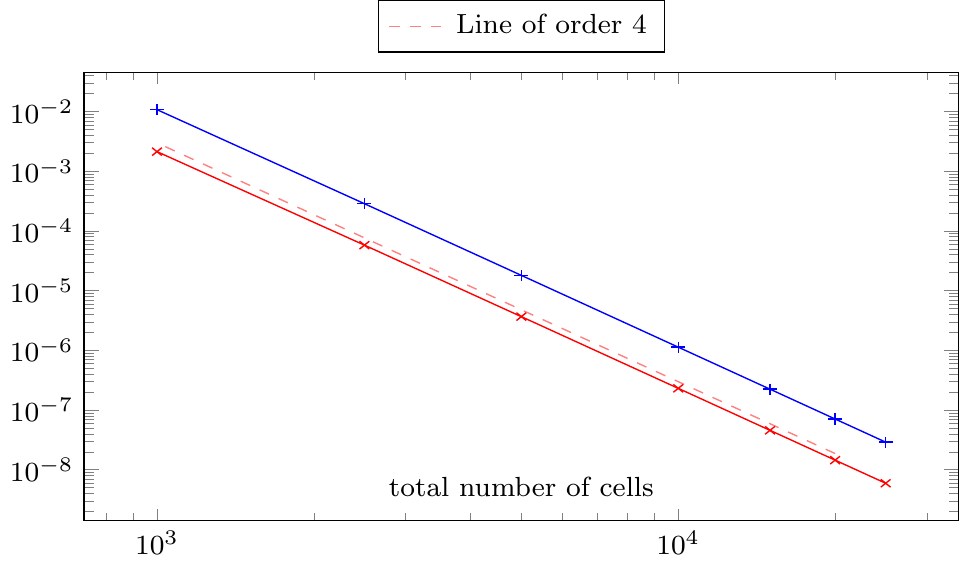}
\caption[a]{Convergence tests for the network problem. The $\ell^1$-error (\includegraphics{tikz/ImplicitActiveFlux-figure_crossref_red.pdf}) and the $\ell^\infty$-error of the point values (\includegraphics{tikz/ImplicitActiveFlux-figure_crossref_blue.pdf}) are shown as a function of the total number of cells (total length of the network is 125). A 4$^\text{th}$ order is shown as dashed reference.
	}
		\label{fig:convergencenetwork}
\end{figure}

Figure \ref{fig:networkconvhalftime} compares the solution on edge 6 at time $t=35$ with results obtained using Methods \exOOLimLLL (4\ref{method:exOOLimLLL}) and \exOLLimLLL (5\ref{method:exOLLimLLL}). At this time the Gaussian, initially present on edge 1 has traveled along edges 3 and $2+5$. On edge 6 one observes the interference between different parts of the Gaussian which have traveled different distances.
After a transitional phase, the steady flow is reached. 
We observe that indeed, the solution on edge 6 is very small.
Figure \ref{fig:networkconv} shows the values on that particular edge for two simulations with $\Delta x = \frac18$ and $\Delta x = \frac1{16}$ using Method \exOLOimLLL (4\ref{method:exOLOimLLL}).
Although inflow at the first edge has amplitude 1, the cancellation is accurate up to $10^{-4}$ and decreasing with resolution. Figure \ref{fig:convergencenetwork} demonstrates experimentally 4$^\text{th}$ order accuracy at $t=100$.

\section{Conclusions and outlook}

Active Flux has been initially introduced in \cite{vanleer77,eymann13} as a fully discrete, time-explicit method, not allowing to be split into a space-discretization and a time-discretization, at least not in any simple way analogous to the method-of-lines. This is possibly the reason for the large domain of stability of this explicit Active Flux method ($\text{CFL} \leq 1$ for 3$^\text{rd}$ order, compare this to DG) and favourable properties concerning the numerical diffusion and dispersion errors (\cite{roe21}). The drawback is that the discretization depends more strongly on the equations to be solved than method-of-lines. 

As part of the quest for even less restrictive CFL-conditions in the present work time-implicit Active Flux methods were obtained. They do not split space and time derivatives and proceed to discretizing them individually, but adopt the genuinely space-time philosophy of the early works \cite{vanleer77,eymann13}. Here, difficulties in following this philosophy were circumvented by means of a reconstruction in time at the interfaces. 
The main contribution of this paper therefore is a rather general way of constructing time-implicit Active Flux methods.

The examples of methods considered in this work are shown to yield compact single-stage methods. They have been shown to perform much better than existing Finite Difference methods of the same orders.
 As commonly observed, even-order methods are found to be more oscillative than odd-order ones, such that the latter are preferred in presence of discontinuities. 
 The ability of the new Active Flux methods to simulate flows on networks has also been demonstrated. 
 The fact that Active Flux has a degree of freedom located at a cell interface makes the implementation of boundary and coupling conditions much easier than for Finite Difference methods. 
 Together, these properties have shown that the new Active Flux methods are very good candidates for industrial applications for network flow.

The time-continuous form for Active Flux (e.g. \cite{abgrall20,abgrall22}) has been only recently obtained and used. Part of the present work was also a comparison between the new methods and time-continuous Active Flux methods endowed with standard time-implicit integration methods. 
The new, genuinely space-time methods presented in this work have been found to perform at least equally well, and many to have significantly less numerical diffusion. 
Due to their single-stage nature, the new Active Flux methods are significantly more efficient than e.g. those using Radau methods, while the DIRK method was much more diffusive than any of the new methods.

Among the methods that were considered in this work, most have a stability domain of $\text{CFL} > 1$, with only one method (4\ref{method:exOOLimLLL}) being unconditionally stable. This method is a good candidate for future development of time-implicit methods for nonlinear problems. Even among the linear ones, not all the possible methods have been studied in the present work, such that in future, more (unconditionally) stable implicit Active Flux methods might be found.

\appendix

\section{Coefficients of the stable methods} \label{app:allmethods}

The methods are ordered by the size of their domain of stability.

\include{allmethods}

\end{document}

%% file: allmethods.tex
    	\subsection{3rd order}
    
    \begin{enumerate}[3A.]
    
    \item \label{method:exOLOimLOL} \exOLOimLOL $\text{CFL} > 1$
    	
\begin{align}
	 0&=  (\cfl-1)  (-1+3  \cfl)	\bar q_i^{n+1}
+ (5+(4-9  \cfl)  \cfl)	\bar q_{i+1}^{n+1} \nonumber \\&
-4q_{i+\frac12}^{n} 	
+(-2+6  \cfl^2) q_{i+\frac32}^{n+1}\\
	0&= (\cfl-1)^2  \cfl^2	\bar q_{i-1}^{n+1}  
 +\cfl^2  (1+\cfl)^2	\bar q_{i+1}^{n+1} \nonumber \\&
 +(2-6  \cfl^2)	\bar q_{i}^{n} 
 -2  (\cfl^2-1) (\cfl^2-1)	\bar q_i^{n+1}\nonumber \\&
-2  \cfl  (\cfl^2-1)q_{i-\frac12}^{n} 	
+2  \cfl  (\cfl^2-1)q_{i+\frac12}^{n} 
\end{align}

    	\item \label{method:exOOLimLOL} \exOOLimLOL $\text{CFL} > 1$
    	 
\begin{align}
0&= 	-4 \bar q_{i+1}^n	
 + \cfl  (-1+3  \cfl)	\bar q_i^{n+1}\nonumber \\&
 - (1+\cfl)  (-4+9  \cfl)	\bar q_{i+1}^{n+1}
+6  \cfl  (1+\cfl) q_{i+\frac32}^{n+1}\\
	0&=  2  (\cfl-1)	\bar q_{i+1}^{n}
 + (\cfl-1)  \cfl	\bar q_{i-1}^{n+1}\nonumber \\&
 -2  (2+\cfl)	\bar q_{i}^{n}
 +(1+\cfl)  (2+\cfl)	\bar q_{i+1}^{n+1} \nonumber \\&
 + (4-2  \cfl  (1+\cfl))\bar q_i^{n+1}
\end{align}

    	\item \label{method:exOLOimOLL} \exOLOimOLL $\text{CFL} > 1$
    	
\begin{align}
	0&=  6  (\cfl-1)  \cfl	\bar q_{i+1}^{n+1} 
+q_{i+\frac12}^n\nonumber \\&
+(-1+(4-3  \cfl)  \cfl)q_{i+\frac12}^{n+1} 	
+(2-3  \cfl)  \cfl q_{i+\frac32}^{n+1}\\
	0&=  - \cfl^3	\bar q_{i+1}^{n+1}
 +(-2+3  \cfl)	\bar q_{i}^{n} \nonumber \\&
 +(2-3  \cfl+\cfl^3)\bar q_i^{n+1} 	
+(\cfl-1)  \cfl q_{i-\frac12}^{n} 	\nonumber \\&
-(\cfl-1)^2  \cfl q_{i-\frac12}^{n+1} 	
-(\cfl-1)  \cfl q_{i+\frac12}^{n} 	\nonumber \\&
+(\cfl-1)^2  \cfl q_{i+\frac12}^{n+1} 
\end{align}

    	\item \label{method:exOOLimOLL} \exOOLimOLL $\text{CFL} > 1$

    	\begin{align}
	0&= \bar q_{i+1}^n	
 + (-1+6  \cfl^2)	\bar q_{i+1}^{n+1}
+(\cfl-3  \cfl^2)q_{i+\frac12}^{n+1} 	\nonumber \\&
-\cfl  (1+3  \cfl) q_{i+\frac32}^{n+1}\\
	0&=  (\cfl-1)  \cfl	\bar q_{i+1}^{n} 
- (1+\cfl)^2	\bar q_{i}^{n} \nonumber \\&
 +\cfl  (1+\cfl)^2	\bar q_{i+1}^{n+1} 
 +(1+2  \cfl-\cfl^2  (2+\cfl))\bar q_i^{n+1} 	\nonumber \\&
+\cfl  (\cfl^2-1)q_{i-\frac12}^{n+1} 	
+(\cfl-\cfl^3)q_{i+\frac12}^{n+1} 
\end{align}

    	\item \label{method:exLOOimLOL} \exLOOimLOL $\text{CFL}\not\in [1, 2]$
    	
\begin{align}
	0&= -4 \bar q_i^n 
+  (13-9  \cfl)  \cfl	\bar q_{i+1}^{n+1}
+ (\cfl-1)  (-4+3  \cfl)	\bar q_i^{n+1} \nonumber \\&
+6  (\cfl-1)  \cfl q_{i+\frac32}^{n+1}\\
	0&=   2  (\cfl-2)	\bar q_{i}^{n}
 +(\cfl-2)  (\cfl-1)\bar q_{i-1}^{n+1} 	\nonumber \\&
 -2  (1+\cfl)	\bar q_{i-1}^{n}
  -2  (\cfl-2)  (1+\cfl)	\bar q_i^{n+1}\nonumber \\&
  +\cfl  (1+\cfl)\bar q_{i+1}^{n+1}
\end{align}

    	\item \label{method:exLOOimOLL} \exLOOimOLL $\text{CFL} > 2$
    	 
\begin{align}
	0&= \bar q_i^n	
 +(5+6  (\cfl-2)  \cfl)	\bar q_{i+1}^{n+1} 
+(-4+(7-3  \cfl)  \cfl)q_{i+\frac12}^{n+1} 	\nonumber \\&
+(2-3  \cfl)  (\cfl-1) q_{i+\frac32}^{n+1}\\
	0&=  (\cfl-2)  (\cfl-1)	\bar q_{i}^{n} 
 - \cfl^2	\bar q_{i-1}^{n}\nonumber \\&
  +(\cfl-1)  \cfl^2	\bar q_{i+1}^{n+1}
 - (\cfl-2)  (-1+\cfl+\cfl^2)	\bar q_i^{n+1}\nonumber \\&
+(\cfl-2)  (\cfl-1)   \cfl q_{i-\frac12}^{n+1} 	
-(\cfl-2)  (\cfl-1)   \cfl q_{i+\frac12}^{n+1} 
\end{align}

    	\item \label{method:exOOLimLLO} \exOOLimLLO $\text{CFL} > 3.74 $
    	 
\begin{align}
	0&= -5 \bar q_{i+1}^n	
  +(-1+6  \cfl^2)	\bar q_i^{n+1}\nonumber \\&
-(1+\cfl)  (-4+9  \cfl)q_{i+\frac12}^{n+1} 	
+(1+\cfl)  (2+3  \cfl) q_{i+\frac32}^{n+1}\\
	0&=  \cfl^2	\bar q_{i+1}^{n} 
  +\cfl^2  (1+\cfl)\bar q_{i-1}^{n+1}	\nonumber \\&
 - (1+\cfl)  (2+\cfl)\bar q_{i}^{n} 	
 +(2-\cfl  (-3+\cfl+\cfl^2))	\bar q_i^{n+1}\nonumber \\&
-\cfl  (1+\cfl)  (2+\cfl)q_{i-\frac12}^{n+1} 	
+\cfl  (1+\cfl)  (2+\cfl)q_{i+\frac12}^{n+1} 
	\end{align}

    	\item \label{method:exOLOimLLO} \exOLOimLLO $\text{CFL} > 4.55$
    	 
\begin{align}
	0&=   6  (\cfl-1)  \cfl	\bar q_i^{n+1}
-5q_{i+\frac12}^{n} 	\nonumber \\&
+(5+(4-9  \cfl)  \cfl)q_{i+\frac12}^{n+1} 	
+\cfl  (2+3  \cfl) q_{i+\frac32}^{n+1} \\
	0&=  (-2-3  \cfl)	\bar q_{i}^{n} 
 + \cfl^3	\bar q_{i-1}^{n+1}\nonumber \\&
+(2+3  \cfl-\cfl^3)	\bar q_i^{n+1}  
-\cfl  (1+\cfl)q_{i-\frac12}^{n} 	\nonumber \\&
-\cfl  (1+\cfl)^2q_{i-\frac12}^{n+1} 	
+\cfl  (1+\cfl)q_{i+\frac12}^{n} 	\nonumber \\&
+\cfl  (1+\cfl)^2q_{i+\frac12}^{n+1} 
	\end{align}

	    	\item \label{method:exLOOimLLO} \exLOOimLLO $\text{CFL} > 4.74$
    	 
    \begin{align}
     0 &= -5 \bar q_i^{n}+ (5+6 (-2+\cfl) \cfl)\bar q_i^{n+1} +(13-9 \cfl) \cfl q_{i+\frac12}^{n+1} \nonumber \\& 
     +\cfl (-1+3 \cfl) q_{i+\frac32}^{n+1}\\
0 &=  (-1+\cfl)^2\bar q_i^{n}+ (-1+\cfl)^2 \cfl \bar q_{i-1}^{n+1} \nonumber \\&
- \cfl (1+\cfl)\bar q_{i-1}^{n} - (1+\cfl) (1+(-3+\cfl) \cfl) \bar q_i^{n+1}\nonumber \\&
+(\cfl-\cfl^3) q_{i-\frac12}^{n+1}+\cfl (-1+\cfl^2) q_{i+\frac12}^{n+1}
    \end{align}

\end{enumerate}

    	\subsection{4th order}
    	
    	\begin{enumerate}[4A.]
    
    \item \label{method:exOOLimLLL} \exOOLimLLL unconditionally stable
    	
    	\begin{align}
	0&= -\bar q_{i+1}^n	
 - \cfl^3	\bar q_i^{n+1}\nonumber \\&
+ (1-\cfl^2  (6+5  \cfl))	\bar q_{i+1}^{n+1} 
+\cfl  (1+\cfl)  (-1+4  \cfl)q_{i+\frac12}^{n+1} 	\nonumber \\&
+\cfl  (1+\cfl)  (1+2  \cfl) q_{i+\frac32}^{n+1}\\
	0&=  2  (\cfl-1)  \cfl	\bar q_{i+1}^{n}
 -2  (1+\cfl)  (2+\cfl)	\bar q_{i}^{n}\nonumber \\&
 + \cfl  (1+\cfl)^2  (2+\cfl)	\bar q_{i+1}^{n+1}
 -2  (\cfl-1)  (2+\cfl)  (1+2  \cfl)	\bar q_i^{n+1}\nonumber \\&
  +(\cfl^2-\cfl^4)	\bar q_{i-1}^{n+1}
+2  (\cfl-1)  \cfl  (1+\cfl)  (2+\cfl)q_{i-\frac12}^{n+1} 	\nonumber \\&
-2  (\cfl-1)  \cfl  (1+\cfl)  (2+\cfl)q_{i+\frac12}^{n+1} 
\end{align}

    	\item \label{method:exOLOimLLL} \exOLOimLLL $\text{CFL} > 1.10$
    	
    	\begin{align}
	0&= \bar q_i^{n+1}  (\cfl-1)  \cfl  (-1+2  \cfl)	
+\bar q_{i+1}^{n+1}  (\cfl-1)  \cfl  (7+10  \cfl)	\nonumber \\&
+2q_{i+\frac12}^{n} 	
-2  (\cfl-1)  (-1+\cfl+4  \cfl^2)q_{i+\frac12}^{n+1} 	\nonumber \\&
+\cfl  (2-4  \cfl^2) q_{i+\frac32}^{n+1}\\
 0&= - (\cfl-1)^2  \cfl^3	\bar q_{i-1}^{n+1}
 +\cfl^3  (1+\cfl)^2	\bar q_{i+1}^{n+1} \nonumber \\&
 +(4-8  \cfl^2)	\bar q_{i}^{n} 
 -4  (\cfl^2-1) (\cfl^2-1)	\bar q_i^{n+1}\nonumber \\&
-2  \cfl  (\cfl^2-1)q_{i-\frac12}^{n} 	
+2  \cfl  (\cfl^2-1) (\cfl^2-1)q_{i-\frac12}^{n+1} 	\nonumber \\&
+2  \cfl  (\cfl^2-1)q_{i+\frac12}^{n} 	
-2  \cfl  (\cfl^2-1) (\cfl^2-1)q_{i+\frac12}^{n+1} 
\end{align}

    	\item \label{method:exLOOimLLL} \exLOOimLLL $\text{CFL} > 1$
    	
\begin{align}
	0&= - \bar q_i^n
 - (\cfl-1)^2 (\cfl-1)	\bar q_i^{n+1}\nonumber \\&
+ \cfl  (-3+(9-5  \cfl)  \cfl)\bar q_{i+1}^{n+1} 	
+(\cfl-1)  \cfl  (-5+4  \cfl)q_{i+\frac12}^{n+1} 	\nonumber \\&
+(\cfl-1)  \cfl  (-1+2  \cfl) q_{i+\frac32}^{n+1}\\
	0&=  2  (\cfl-2)  (\cfl-1)	\bar q_{i}^{n} 
 - (\cfl-2)  (\cfl-1)^2  \cfl	\bar q_{i-1}^{n+1}\nonumber \\&
 -2  \cfl  (1+\cfl)	\bar q_{i-1}^{n}
 -2  (\cfl-2)  (1+\cfl)  (-1+2  \cfl)	\bar q_i^{n+1}\nonumber \\&
 + \cfl^2  (\cfl^2-1)	\bar q_{i+1}^{n+1}
+2  (\cfl-2)  (\cfl-1)  \cfl  (1+\cfl)q_{i-\frac12}^{n+1} 	\nonumber \\&
-2  (\cfl-2)  (\cfl-1)  \cfl  (1+\cfl)q_{i+\frac12}^{n+1} 
\end{align}
		
	\item \label{method:exLOLimLOL} \exLOLimLOL marginally stable
    	
\begin{align}
	0&=  (\cfl-1)  (-5+4  \cfl)	\bar q_{i+1}^{n} 
 - (1+\cfl)  (-1+4  \cfl)	\bar q_{i}^{n}\nonumber \\&
 +(\cfl-1)  (1+\cfl  (-5+3  \cfl))	\bar q_i^{n+1} 
 - (1+\cfl)  (5+\cfl  (-17+9  \cfl))	\bar q_{i+1}^{n+1}\nonumber \\&
+6  \cfl  (\cfl^2-1) q_{i+\frac32}^{n+1}\\
	 0&=  (\cfl-2)  (\cfl-1)	\bar q_{i-1}^{n+1}
 - (\cfl-2)  (\cfl-1)	\bar q_{i+1}^{n}\nonumber \\&
 - (1+\cfl)  (2+\cfl)	\bar q_{i-1}^{n}
 + (1+\cfl)  (2+\cfl)	\bar q_{i+1}^{n+1}\nonumber \\&
 +(8-2  \cfl^2)	\bar q_i^{n+1} 
 +2  (-4+\cfl^2)\bar q_{i}^{n}
\end{align}
    	
    	\end{enumerate}
    	
    	       \subsection{5th order}

    	\begin{enumerate}[5A.]
    
    \item \label{method:exOLLimLLL} \exOLLimLLL $\text{CFL} > 1$

    	\begin{align}
	 0&= - (\cfl-1)  \cfl^2  (-1+5  \cfl^2) \bar q_i^{n+1}	
 + 2 (\cfl-1)  (-2+\cfl+5  \cfl^2)	\bar q_{i+1}^{n} \nonumber \\&
  -(\cfl-1)  (1+\cfl)^2  (-4+5  \cfl  (2+5  \cfl))\bar q_{i+1}^{n+1}	
-2  \cfl  (1+\cfl)  (1+5  \cfl)q_{i+\frac12}^{n} 	\nonumber \\&
+2  (\cfl-1)  \cfl  (1+\cfl)  (-3+5  \cfl  (1+2  \cfl))q_{i+\frac12}^{n+1} 	
+2  \cfl  (1+\cfl)^2  (-2+5  \cfl^2) q_{i+\frac32}^{n+1}  \\
	  0&=  (\cfl-1)^2  \cfl^2	\bar q_{i+1}^{n}
 + \frac12   (\cfl-1)^2  \cfl^3	\bar q_{i-1}^{n+1} \nonumber \\&
  -\frac12   \cfl^2  (1+\cfl)^2  (2+\cfl)	\bar q_{i+1}^{n+1}
  +(\cfl-1)^2  (2+\cfl)  (2+3  \cfl)	\bar q_i^{n+1}\nonumber \\&
 +(-4+\cfl^2  (9-(\cfl-2)  \cfl))	\bar q_{i}^{n} 
+(\cfl-1)  \cfl  (1+\cfl)  (2+\cfl)q_{i-\frac12}^{n} 	\nonumber \\&
-(\cfl-1)^2  \cfl  (1+\cfl)  (2+\cfl)q_{i-\frac12}^{n+1} 	
-(\cfl-1)  \cfl  (1+\cfl)  (2+\cfl)q_{i+\frac12}^{n} 	\nonumber \\&
+(\cfl-1)^2  \cfl  (1+\cfl)  (2+\cfl)q_{i+\frac12}^{n+1} 
\end{align}

    	\item \label{method:exLLOimLLL} \exLLOimLLL  $\text{CFL} > 2 $

    	\begin{align}
	0&=  - (\cfl-1)^2  (4+5  (\cfl-2)  \cfl) 	\bar q_i^{n+1}
  -2  (-2+\cfl+5  \cfl^2) \bar q_{i}^{n} 	\nonumber \\&
 +\cfl  (-16+\cfl  (9+5  (8-5  \cfl)  \cfl)) 	\bar q_{i+1}^{n+1} 
+2  (\cfl-1)  (-4+5  \cfl)q_{i+\frac12}^{n}  	\nonumber \\&
+2  (\cfl-1)  (-1+2  \cfl)  (-4+5  (\cfl-1)  \cfl)q_{i+\frac12}^{n+1}  	
+2  (\cfl-1)  \cfl  (-2+5  \cfl^2) q_{i+\frac32}^{n+1}\\
	0&=  (\cfl-2)  (\cfl-1)^2  \cfl^2 	\bar q_{i-1}^{n+1} 
  +2  \cfl^2  (1+\cfl)^2 \bar q_{i-1}^{n} 	\nonumber \\&
 - \cfl^3  (1+\cfl)^2 	\bar q_{i+1}^{n+1}
 +2  (\cfl-2)  (1+\cfl)^2  (-2+3  \cfl) \bar q_i^{n+1} 	\nonumber \\&
 -2  (4+\cfl^2  (-9+\cfl  (2+\cfl))) \bar q_{i}^{n} 	
-2  (\cfl-2)  (\cfl-1)  \cfl  (1+\cfl)q_{i-\frac12}^{n}  	\nonumber \\&
-2  (\cfl-2)  (\cfl-1)  \cfl  (1+\cfl)^2q_{i-\frac12}^{n+1}  	
+2  (\cfl-2)  (\cfl-1)  \cfl  (1+\cfl)q_{i+\frac12}^{n}  	\nonumber \\&
+2  (\cfl-2)  (\cfl-1)  \cfl  (1+\cfl)^2q_{i+\frac12}^{n+1}
\end{align}

    	\item \label{method:exLOLimLLL} \exLOLimLLL $\text{CFL} > 2$
    	
  \begin{align}  	
	0&=  (\cfl-1)^2  (-4+5  \cfl)	\bar q_{i+1}^{n} 
 - \cfl  (1+\cfl)  (1+5  \cfl)\bar q_{i}^{n}	\nonumber \\&
  -(\cfl-1)^2  \cfl  (-1+5  (\cfl-1)  \cfl)	\bar q_i^{n+1}
 - (1+\cfl)  (-4+\cfl  (17+\cfl  (-4+5  \cfl  (-8+5  \cfl))))	\bar q_{i+1}^{n+1}\nonumber \\&
+2  (\cfl-1)  \cfl  (1+\cfl)  (2+5  \cfl  (-3+2  \cfl))q_{i+\frac12}^{n+1} 	
+2  \cfl  (2-7  \cfl^2+5  \cfl^4) q_{i+\frac32}^{n+1}\\
	0&=  - (\cfl-2)  (\cfl-1)^2  \cfl	\bar q_{i+1}^{n}
 - (\cfl-2)  (\cfl-1)^2  \cfl  (1+\cfl)	\bar q_{i-1}^{n+1}\nonumber \\&
 - \cfl  (1+\cfl)^2  (2+\cfl)	\bar q_{i-1}^{n}
+ (\cfl-1)  \cfl  (1+\cfl)^2  (2+\cfl)	\bar q_{i+1}^{n+1} \nonumber \\&
+ (-8+26  \cfl^2-6  \cfl^4)	\bar q_i^{n+1} 
+ 2  (4-5  \cfl^2+\cfl^4)	\bar q_{i}^{n}\nonumber \\&
+2  \cfl  (4-5  \cfl^2+\cfl^4)q_{i-\frac12}^{n+1} 	
-2  \cfl  (4-5  \cfl^2+\cfl^4)q_{i+\frac12}^{n+1} 
\end{align}

    	\end{enumerate}